\documentstyle[10pt,intlim,righttag]{amsart}

\newcommand{\essinf}{\operatornamewithlimits{ess\,inf}}
\newcommand{\esssup}{\operatorname{ess\,sup}}
\newcommand{\diag}{\operatorname{diag}}

\newcommand{\vp}{\varphi}
\newcommand{\vps}{\varepsilon}
\newcommand{\td}{\tilde}

\newcommand{\lbl}{\label}

\newcommand{\la}{\langle}
\newcommand{\ra}{\rangle}

\newcommand{\wt}{\widetilde}

\newcommand{\be}{\begin}
\newcommand{\ee}{\end}

\newtheorem{prop}{\ \ \ Proposition}[section]

\theoremstyle{definition}

\renewcommand{\rm}{\normalshape}

\numberwithin{equation}{section}

\begin{document}

\title[Backward Stochastic PDE]
{Backward Stochastic PDE$s$ related to the \\utility maximization
problem}

\author[Michael Mania and Revaz Tevzadze]{Michael Mania$^{1),2)}$ and Revaz Tevzadze$^{1),3)}$}
\maketitle

\begin{center}
$^{1)}$ Georgian--American University, Business School,
3, Alleyway II, \\[2mm]  Chavchavadze Ave. 17\,a, Tbilisi, Georgia,\\[2mm]
$^{2)}$ A. Razmadze Mathematical Institute, 1, M. Aleksidze St.,
Tbilisi, Georgia, E-mail:mania@@rmi.acnet.ge \\[2mm]
$^{3)}$Institute of Cybernetics, 5, S. Euli St., Tbilisi, Georgia,
E-mail:tevza@@cybernet.ge
\end{center}

\begin{abstract}

We study  utility maximization problem for general utility functions
using dynamic programming approach. We consider an incomplete
financial market model, where the dynamics of asset prices are
described by an $R^d$-valued continuous semimartingale. Under some
regularity assumptions we derive backward stochastic partial differential equation
(BSPDE) related
directly to the primal problem and show  that the strategy is optimal
if and only if the corresponding wealth process satisfies a certain forward-SDE.
 As examples  the cases of power,
exponential and logarithmic utilities
are considered.

\bigskip

\noindent {\bf 2000  Mathematics Subject Classification}:
90A09,60H30, 90C39.

\noindent {\bf Key words and phrases}: Backward stochastic partial
differential equation, utility maximization problem, semimartingale, incomplete markets.

\end{abstract}

\section{Introduction}

\

Portfolio optimization, hedging and derivative pricing are
fundamental problems in mathematical finance, which are closely
related to each other. A basic optimization problem of
mathematical finance, such as optimal portfolio choice or hedging,
is to optimize
\begin{equation}\label{ut1}
 E[ U(X^{x,\pi}_T)]\;\;\;\;\text{over
all}\;\;\;\;\pi\;\;\;\text{from a class}\;\;\Pi\;\;\;\text {of
strategies},
\end{equation}
where $X^{x,\pi}_t=x+\int_0^t\pi_udS_u$ is the wealth process
starting from initial capital $x$, determined by the
self-financing trading strategy $\pi$ and $\Pi$ is some class of
admissible strategies. $U$ is an objective function which can be
depended also on $\omega$. It can be interpreted as a utility
function or a function which measures a hedging error.

If we assume that $U(x)$ is  strictly convex  (for each $\omega$)
then one can  interpret $U$ as a function which measures a hedging
error and consider the  problem
\begin{equation}\label{utt1}
\text {to minimize}\;\;\; E[ U(X^{x,\pi}_T)]\;\;\;\;\text{over
all}\;\;\;\;\pi\;\;\;\text{from }\;\;\Pi.
\end{equation}
In \cite{MT7}
a backward stochastic PDE for value function
\begin{equation}\label{utt2}
V(t,x)=\underset{\pi\in\Pi}{\text{\rm essinf}}\text{ }
E(U(x+\int_t^T\pi_udS_u )/{\mathcal F}_t)
\end{equation}
of (\ref{utt1}) was derived and in terms of solutions of this equation
a characterization of optimal strategies was given. We shall use the same
approach to the case when the objective function $U$ is strictly concave. In particular,
if $U$ is a utility function,
then (\ref{ut1}) corresponds to the utility  maximization problem
 \begin{equation}\label{ut3}
\text{to maximize}\;\;\;\;\; E[ U(X^{x,\pi}_T)]\;\;\;\;\text{over
all}\;\;\;\;\pi\in\Pi,
\end{equation}
i.e., for a given initial capital $x>0$, the goal is to maximize
the expected value from the terminal wealth.

The utility maximization problem was first studied by Merton
(1971) in the classical Black-Scholes model. Using the Markov
structure of the model he derived the Bellman equation for the
value function of the problem and produced the closed-form
solution of this equation in cases of power, logarithmic and
exponential utility functions.

For general complete market models, it was shown by Pliska (1986),
Cox and Huang (1989) and Karatzas et al (1987) that the optimal
portfolio of the utility maximization problem is (up to a
constant) equal to the density of the martingale measure, which is
unique for complete markets. As shown by He and Pearson (1991) and
Karatzas et al (1991), for incomplete markets described by
Ito-processes, this method gives a duality characterization of
optimal portfolios provided by the set of martingale measures.
Their idea was to solve the dual problem of finding the suitable
optimal martingale measure and then to express the solution of the
primal problem by convex duality. Extending the domain of the dual
problem the approach has been generalized to semimartingale models
and under weaker conditions on the utility functions by Kramkov
and Schachermayer (1999). See also more recent papers \cite{Fr},
\cite{El-R}, \cite{CSW}, \cite{Mnif},
\cite{GR 01}, \cite{D-Gr-Rh-S-S-S}, \cite{Sch}, \cite{k-sc2}, \cite{BF}.

These approaches mainly give a reduction of the basic primal
problem to the solution of the  dual problem, but the constructive
solution of the dual problem for general models of incomplete
markets is itself demanding task.

Our goal is to derive a semimartingale Bellman equation (a
stochastic version of the Bellman
 equation) related directly to the basic (or primal) optimization problem
 and to give constructions of
 optimal strategies. Applying the dynamic programming approach directly
 to the primal optimization problem may in many cases represent a valuable
 alternative to the commonly used convex duality approach.

Let $S$ be an $R^d$-valued continuous semimartingale, defined on a
filtered probability space satisfying the usual conditions. The
process $S$ describes the discounted price evolution of $d$ risky
assets in a financial market containing also a riskless bond with
a constant price. To exclude arbitrage opportunities, we suppose
that the set ${\mathcal M}^e$ of equivalent martingale measures
for $S$ is not empty. Since $S$ is continuous, the existence of an
equivalent martingale measure implies that the structure condition
is satisfied, i.e., $S$ admits the decomposition
\begin{equation}\label{ut4}
S_t=M_t+\int_0^td\langle M\rangle_s\lambda_s,\;\; \int_0^t
\lambda'_sd\la M\ra_s\lambda_s <\infty\;\;
\text{for all} \;\; t\;\;a.s.,
\end{equation}
where $M$ is a continuous local martingale and $\lambda$ is a
predictable $R^d$-valued process.

We consider utility function $U$ mapping $R_+\equiv(0,\infty)$ into $R$. It
is assumed to be continuously differentiable, strictly increasing,
strictly concave and to satisfy the Inada conditions:
$$
U^{'}(0)=\lim_{x\to 0}U^{'}(x)=\infty,
$$
$$
U^{'}(\infty)=\lim_{x\to\infty}U^{'}(x)=0.
$$
We also set $U(0)=\lim_{x\to0}U(x)$ and $U(x)=-\infty$ for all
$x<0$.

Denote by ${\cal M}^e$ the set of martingale measures for $S$. Throughout the paper
we assume that
$$
{\cal M}^e\neq\emptyset.
$$

For any $x\in R_+$, we denote by $\Pi_x$ the class of predictable
$S$-integrable processes $\pi$ such that the corresponding wealth
process is nonnegative  at any instant, that is
$X^{x,\pi}_t=x+\int_0^t\pi_udS_u\ge0$ for all $t\in [0,T]$.

For simplicity in introduction we consider the case with one risky
asset.

Let us introduce the dynamical value function of the problem
(\ref{ut3}) defined as
\begin{equation}\label{ut5}
V(t,x)=\underset{\pi\in\Pi_x}{\text{\rm esssup}}\text{ }
E\big(U(x+\int_t^T\pi_udS_u )/{\mathcal F}_t\big).
\end{equation}

The classical It\^o formula (or its generalization by Krylov 1980)
plays a crucial role to derivation of the Bellman equation for the
value function of controlled diffusion processes. For our purposes
the It\^o formula is not sufficient since the function $V$ depends
also on $\omega$, even if $U$ is deterministic. Therefore the
It\^o-Ventzel formula should be used.

Under some regularity assumptions on the value function
(sufficient for the application of the It\^o--Ventzell formula) we
show in Theorem \ref{utt3.1} that value function defined by
(\ref{ut5}) satisfies the following backward stochastic partial
differential equation (BSPDE)
\begin{equation}\label{ut6}
V(t,x)= V(0,x)+\frac{1}{2}\int_0^t
\frac{(\varphi_x(s,x)+\lambda(s)V_x(s,x))^2} {V_{xx}(s,x)}d\langle M\rangle_s
\end{equation}
$$
+\int_0^t\varphi(s,x)dM_s+L(t,x)
$$
with the boundary condition
$$
V(T,x)=U(x),
$$
where $\int_0^t\varphi(s,x)dM_s+L(t,x)$ is the  martingale part of $V(t,x)$, $L(t,x)$ is
strongly orthogonal to $M$ for all $x$
and subscripts $\varphi_x, V_x, V_{xx}$ stand for the partial
derivatives. Moreover, the strategy $\pi^*$ is optimal if and only
if the corresponding wealth process $X^{\pi^*}$ is a solution of
the following forward SDE
\begin{equation}\label{ut7}
X_t^{\pi^*}=X_0^{\pi^*}- \int_0^t\frac{\varphi_x(u,X_u^{\pi^*})+
\lambda(u)V_x(u, X_u^{\pi^*})} {V_{xx}(s,X_u^{\pi^*})} dS_u.
\end{equation}

Thus, to give the construction of the optimal strategy one
should:
\\ \noindent1) first solve the backward equation
(\ref{ut6}) (which determines $V$ and $\varphi$ simultaneously)
and substitute corresponding derivatives of $V$ and $\varphi$ in
equation (\ref{ut7}), then
\\ \noindent
2) solve the forward equation (\ref{ut7}) with respect to
$X^{\pi^*}$ and, finally,
\\ \noindent
3) reproduce the optimal strategy $\pi^*$ from the corresponding
wealth process $X^{\pi^*}$.

Theorem \ref{utt3.1} is a verification theorem, since we require
conditions directly on the value function $V$ and not only on the
basic objects (on the model and on the objective function $U$). Therefore we can't
state that the solution of
equation (\ref{ut6}) exists, but for standard utility functions
(e.g., for power, exponential, logarithic and quadratic utilities)
all conditions of Theorem \ref{utt3.1} are satisfied and in these
cases the existence of a unique solutions of corresponding
backward equations follows from this theorem.

 If $U(x)=x^p, p\in (0,1)$,
then (\ref{ut3}) corresponds to power utility maximization problem
\begin{equation}\label{ut8}
\text{to maximize}\;\;\;\;\;
E(x+\int_0^T\pi_udS_u)^p\;\;\;\;\text{over
all}\;\;\;\;\pi\in\Pi_x.
\end{equation}

In this case $V(t,x)=x^pV_t$, where $V_t$ is a semimartingale and
all condition of Theorem 3.1 are satisfied. This theorem implies
that the process $V_t$ satisfies the following backward stochastic
differential equation (BSDE)
$$
V_t= V_0+\frac{q}{2}\int_0^t \frac{(\varphi_s+\lambda_sV_s)^2}
{V_s}d\langle M\rangle_s
$$
\begin{equation}\label{ut9}
+\int_0^t\varphi_sdM_s+L_t,\;\;V_T=1,
\end{equation}
\\ \noindent
where $q=\frac{p}{p-1}$ and $L$ is a local martingale strongly
orthogonal to $M$.
\\ \noindent
Besides, equation (\ref{ut7}) is transformed into a linear
equation
\begin{equation}\label{ut10}
X_t^{*}=x+(1-q)\int_0^t\frac{\varphi_u+\lambda_uV_u}
{V_u}X^{*}_udS_u
\end{equation}
for the optimal wealth process.

Therefore,
$$
X_t^{*}=x{\mathcal E}_t((1-q)(\frac{\varphi}{V}+\lambda)\cdot S)
$$
and the optimal strategy is of the form
$$
\pi_t^*=x(1-q)(\frac{\varphi_t}{V_t}+\lambda_t){\mathcal
E}_t\big((1-q)(\frac{\varphi}{V}+\lambda)\cdot S\big).
$$
Equations of type (\ref{ut9}) was derived in \cite{MRT} in relation
to utility maximization problem and in \cite{Hu} for constrained utility
maximization problem. In comparison to the work \cite{Hu} our results are
at the same time more and less general. In \cite{Hu} diffusion market model
is considered and the boundedness of model coefficients is assumed.
 We are working  with a general right-continuous filtration and under weaker
 boundedness conditions, but we have not included constraints on our strategies.

We consider also utility functions which take finite values on all real line, such as
the exponential utility $U(x)=1-e^{-\gamma x}$.  For this case Theorem 3.1 is not
directly applicable. It needs a special choice of a class of trading strategies
 and additional assumption of the existence of the dual optimizer (see Schachermayer 2003).
Exponential utility maximization problem  we consider in section 4 and
distinguish cases when this problem admits an explicit solution. In section 4 we consider
also the case of quadratic utility.

The main tools of the work - Backward Stochastic Differential
Equations, have been introduced by J. M. Bismut in \cite{Bs}  for
the linear case as the equations for the adjoint process in the
stochastic maximum principle. In \cite{Ch} and \cite{Par-P}
the well-posedness results for BSDEs with more general generators
was obtained (see also \cite{El-P-Q} for references and related
results).  The semimartingale backward equation, as a stochastic
version of the Bellman equation in an optimal control problem, was
first derived in \cite{Ch} by R. Chitashvili.

The main results of this paper are based on the papers of authors
 \cite{MT7}, \cite{MRT}.

\

\section{Basic assumptions and some auxiliary facts}

\

We consider an incomplete financial market model, where the
dynamics of asset prices are described by an $R^d$-valued
continuous semimartingale $S$  defined on a filtered probability
space $(\Omega,F,{\cal{F}}=({\cal F}_t,t\in[0,T]),P)$ satisfying
the usual conditions, where $F={\cal F}_T$ and $T<\infty$ is a
fixed time horizon.
For all unexplained notations concerning the martingale theory used below we
refer the reader to \cite{J},\cite{D-M},\cite{L-Sh2}.

Denote by ${\mathcal M}^e$  the set of martingale measures, i.e.,
the set of measures $Q$ equivalent to $P$ on ${\mathcal F}_T$ such
that $S$ is a local martingale under $Q$. Let $Z_t(Q)$ be the
density process of $Q$ with respect to the basic measure $P$,
which is a strictly positive uniformly integrable martingale. For
any $Q\in{\cal M}^e$ there is a $P$-local martingale $M^Q$ such
that $Z(Q)={\cal{E}}(M^Q)=({\cal{E}}_t(M^Q),t\in[0,T]),$ where
${\cal{E}}(M)$ is the Doleans-Dade exponential of $M$.

We recall the definition of {\rm BMO}-martingales and the Muckenhoupt
condition.

The square integrable continuous martingale $M$ belongs to the
class {\rm BMO} if there is a constant $C>0$ such that $$
E(\langle M\rangle_T-\langle M\rangle_\tau|F_{\tau})\le C, \:\:\:\:\:P-{\mbox a.s.} $$
for every stopping time $\tau$.

A strictly positive uniformly integrable martingale $Z$ satisfies
the Muckenhoupt inequality denoted by $A_\alpha(P)$ for some $1<\alpha<\infty$,
iff there is a constant $C$ such that
$$
E(\big(\frac{Z_\tau}{Z_T}\big)^{\frac{1}{\alpha-1}}|F_{\tau})\le C,
\:\:\:\:\:P-{\mbox a.s.}
$$
for every stopping time $\tau$.

Note that, if the mean variance tradeof $\la \lambda\cdot M\ra_T$ is bounded, then
the density process ${\cal E}(-\lambda\cdot M)$ of the minimal martingale measure
satisfies   the Muckenhoupt inequality for any $\alpha>1$.

The following assertion relates $BMO$ and the Muckenhoupt
condition.

 \begin{prop}(\cite{Dol-M}, \cite{Kz}). Let $M$ be a local
 martingale and $\mathcal E(M)$ its Dol\'eans Exponential. The following
 assertions are equivalent:

 (i) $M$ belongs to the class $BMO$,

 (ii) $\mathcal E(M)$ is a uniformly integrable martingale satisfying
 the Muckenhoupt inequality
 $A_\alpha(P)$ for some $\alpha>1$.
 \end{prop}

Let $\Pi_x$ be the space of all predictable $S$-integrable
processes $\pi$ such that the corresponding wealth process is
nonnegative at any instant, that is $x+\int_0^t\pi_udS_u\ge0$ for
all $t\in[0,T]$.

In the sequel  sometimes we shall use the notation $(\pi\cdot
S)_t$ for the stochastic integral $\int_0^t\pi_udS_u$.

By $\mu^\kappa$ we shall denote the Dolean's measure of an increasing process $\kappa$.

Suppose that the objective function $U(x)=U(\omega,x)$ satisfies
the following conditions:

{\bf B1)}  $V(0, x)<\infty$ for some $x$,

{\bf B2)} $U(\omega, x)$ is  twice continuously differentiable and strictly concave
for each $\omega$,

{\bf B3)} optimization problem (\ref{ut3}) admits a solution,
i.e., for any $t$ and $x$ there is a strategy  $\pi^*(t,x)$ such
that
\begin{equation}\label{ut13}
V(t,x)=E(U(x+\int_t^T\pi^*_s(t,x)dS_s )/{\cal F}_t)
\end{equation}

{\bf Remark 2.1}. As shown by Kramkov and Schachermayer (1999),
sufficient  condition for B3), when $U$ does not depend on $\omega$, is that utility
function $U(x)$ has asymptotic elasticity strictly less than 1,
i.e.,
\begin{equation}\label{ae}
AE(U)=\limsup_{x\to\infty}\frac{xU_x(x)}{U(x)}<1.
\end{equation}
It follows from Kramkov and Schachermayer (2003) that for B3)  the
finitness of the dual value function is also sufficient.

{\bf Remark 2.2}. The strict concavity of $U$ implies that the
optimal strategy is unique if it exists. Indeed, if there exist
two optimal strategies $\pi^1$ and $\pi^2$, then by concavity of
$U$ the strategy $\bar\pi=\frac{1}{2}\pi^1+\frac{1}{2}\pi^2$ is
also optimal. Therefore,
$$
\frac{1}{2}E[U(x+\int_t^T{\pi^1_s}dS_s )|{\cal F}_t]+\frac{1}{2}
E[U(x+\int_t^T{\pi^2_s}dS_s )|{\cal F}_t]
$$
$$
= E[U(x+\int_t^T{\bar\pi_s}dS_s )|{\cal F}_t]
$$
and
$$
\frac{1}{2}U(x+\int_t^T{\pi^1_s}dS_s )+
\frac{1}{2}U(x+\int_t^T{\pi^2_s}dS_s )
$$
$$
 =U(x+\int_t^T{\bar\pi_s}dS_s )
\;\;P-a.s.
$$
Now strict concavity of $U$ leads to the equality
$\int_t^T{\pi^1_s}dS_s = \int_t^T{\pi^2_s}dS_s$.

  For convenience we give the
proof of the following known assertion.

\be{lem}\lbl{utl2.1} Under conditions B1)-B3) the value function
$V(t,x)$ is a strictly concave function with respect to $x$.
\ee{lem}

{\it Proof}. The  concavity of $V(t,x)$ follows from B2) and B3),
since for any $\alpha,\beta\in [0,1]$ with $\alpha+\beta=1$ and
any $x_1,x_2\in R$ we have
$$
\alpha V(t,x_1)+\beta V(t,x_2)
$$
$$
=\alpha E[U(x_1+\int_t^T\pi^*_u(t,x_1)dS_u)|{\cal F}_t]+ \beta
E[U(x_2+\int_t^T\pi^*_u(t,x_2)dS_u)|{\cal F}_t]
$$
$$
\ge E[U(\alpha x_1+\beta x_2+\int_t^T(\alpha\pi_u^*(t,x_1)+
\beta\pi^*_u(t,x_2))dS_u|{\cal F}_t]
$$
\begin{equation}\label{ut14}
\ge V(t,\alpha x_1+\beta x_2).
\end{equation}
To show that $V(t,x)$ is strictly concave we must verify that if
the equality
\begin{equation}\label{ut15}
\alpha V(t,x_1)+\beta V(t,x_2)= V(t,\alpha x_1+\beta x_2)
\end{equation}
is valid for some $\alpha,\beta\in (0,1)$ with $\alpha+\beta=1$,
then $x_1=x_2$.

Indeed, if equality (\ref{ut15}) holds, then from (\ref{ut14}) and
the strict convexity of $U$ follows that $P$-a.s.
$$
x_1+\int_t^T\pi^*_u(t,x_1)dS_u= x_2+\int_t^T\pi^*_u(t,x_2)dS_u,
$$
which implies that $x_1=x_2$.

{\bf Remark 2.3.} The concavity of $V(0,x)$ and condition $B1)$
imply that $V(0,x)<\infty$ for all $x\in R_+$.

{\it Ito-Ventzell's formula.}

Let $(Y(t,x),\;\;t\in[0,T],\;\;x\in R)$ be a family of special
semimartingales with the decomposition
\begin{equation}\label{ut16}
Y(t,x)=Y(0,x)+B(t,x)+N(t,x),
\end{equation}
where $B(\cdot,x)\in{\cal A}_{\text loc}$ and $N(\cdot,x)\in{\cal
M}_{\text{loc}}$ for any $x\in R$. By the
Galtchouk--Kunita--Watanabe (G-K-W) decomposition of $N(\cdot,x)$
with respect to $M$  a parametrized family of semimartingales $Y$
admits the representation
\begin{equation}\label{ut17}
Y(t,x)=Y(0,x)+B(t,x)+ \int_0^t\psi(s,x)dM_s+L(t,x),
\end{equation}
where $L(\cdot,x)$ is a local martingale strongly orthogonal to
$M$ for all $x\in R$.

Assume that:

C1) there exists a predictable increasing process
$(K_t,t\in[0,T])$ such that $B(\cdot,x)$ and $\langle M\rangle$ are absolutely
continuous with respect to $K$, i.e., there is a measurable
function $b(t,x)$ predictable for every $x$ and a matrix-valued
predictable process $\nu_t$ such that
$$
B(t,x)=\int_0^tb(s,x)dK_s,\;\;\langle M\rangle_t=\int_0^t\nu_sdK_s.
$$
Note that, by continuity of $M$ the square characteristic $\langle M\rangle$ is
absolutely continuous w.r.t. the continuous part $K^c$ of the
process $K$ and
$$
\langle M\rangle_t=\int_0^t\nu_sdK^c_s=\int_0^t\nu_sdK_s.
$$
Without loss of generality one can assume that $\nu$ is bounded
and the scalar product $u'\nu_tv$ for $u,v\in R^d$ we denote by
$(u,v)_{\nu_t}$ .

C2) the mapping $x\to Y(t,x)$ is twice continuously differentiable
for all $(\omega,t)$,

C3) the first derivative $Y_x(t,x)$ is a special semimartingale,
admitting
 the decomposition
\begin{equation}\label{ut18}
Y_x(t,x)=Y_x(0,x)+B_{(x)}(t,x)+\int_0^t\psi_x(s,x)dM_s+L_{(x)}(t,x),
\end{equation}
where $B_{(x)}(\cdot,x)\in{\cal A}_{\text loc}$,
$L_{(x)}(\cdot,x)$ is a local martingale orthogonal to $M$ for all
$x\in R$ and $\psi_x$ is the partial derivative of $\psi$ at $x$
(note that $A_{(x)}$ and $L_{(x)}$ are not assumed to be
derivatives of $A$ and $L$ respectively, whose existence does not
necessarily follow from condition C2)),

C4) $Y_{xx}(t,x)$ is RCLL process for every $x\in R$,

C5) the functions $b(s,\cdot),\;\;\psi(s,\cdot)$ and
$\psi_x(s,\cdot)$ are continuous at $x$  $\mu^{K}$-a.e.,

C6) for any $c>0$
$$
E\int_0^T\sup_{|x|\le c} g(s,x)dK_s<\infty
$$
for $g$ equal to $|b|,\;|\psi|^2$ and $|\psi|^2_x$.


In what follows we shall need the following version of
Ito-Ventzell's formula

\be{prop}\lbl{utp2.1} Let $(Y(\cdot,x),\;\;x\in R)$ be a family of
special semimartingales satisfying conditions C1)-C6) and
$X^\pi=x+\pi\cdot S$. Then the transformed process
$Y(t,X_t^\pi),\;t\in[0,T]$ is a special semimartingale with the
decomposition
$$
 Y(t,X_t^\pi)=Y(0,c)+B_t+N_t,
$$
where
$$
B_t=\int_0^t\big[Y_x(s,X_s^\pi)\lambda_s'd\langle M\rangle_s\pi_s+
\psi_x(s,X_s^\pi)'d\langle M\rangle_s\pi_s+
$$
\begin{equation}\label{ut19}
+\frac{1}{2}Y_{xx}(s,X_s^\pi)\pi_s'd\langle M\rangle_s\pi_s\big]+\int_0^tb(s,X_s^\pi)dK_s
\end{equation}
and $N$ is a continuous local martingale. \ee{prop}

\noindent
One can derive this assertion from Theorem 1.1 of
\cite{L-T} or from Theorem 2 of \cite{Ch-M}. Here we don't require
any conditions on $L(t,x)$ imposed in \cite{L-T} and \cite{Ch-M},
since the martingale part of substituted process $X^\pi$ is
orthogonal to $L(\cdot,x)$ and since we don't give an explicit
expression of martingale part $N$, which is not necessary for our
purposes.

{\bf Remark 2.4.} Since the semimartingale $S$ is assumed to be
continuous and is of the form (\ref{ut4}), only the latter term of
(\ref{ut19}) may have the jumps, i.e., the process $K$ is not
continuous in general.

\

\section{The BSPDE for the value function}

\

In this section we derive a backward stochastic PDE for the value
function related to the utility maximization problem.

Denote by ${\mathcal V}^{1,2}$ the class of functions
$Y:\Omega\times[0,T]\times R\to R$ satisfying conditions C1)-C6).

Let us consider the following backward stochastic partial
differential equation (BSPDE)
\begin{align}
\notag Y(t,x)&= Y(0,x)
\\
\notag
&+\frac{1}{2}\int_0^t \frac{(\psi_x(s,x)+\lambda(s)Y_x(s,x))'}
{Y_{xx}(s,x)}d\langle M\rangle_s(\psi_x(s,x)+\lambda(s)Y_x(s,x))
\\
\label{ut20}
&+\int_0^t\psi(s,x)dM_s+L(t,x),\;\;\;\;\;\;\;\;\;\;\;\;\;\;\;\; L(\cdot,x)\bot M,
\end{align}

with the boundary condition
\begin{equation}\label{ut21}
Y(T,x)=U(x).
\end{equation}

We shall say that $Y$ solves equation (\ref{ut20}),(\ref{ut21})
if:

(i) $Y(\omega,t,x)$ is twice continuously differentiable for each
$(\omega,t)$ and satisfies the boundary condition (\ref{ut21}),

(ii) $Y(t,x)$ and $Y_x(t,x)$ are special semimartingales admitting
decompositions (\ref{ut17}) and (\ref{ut18}) respectively, where
$\psi_x$ is the partial derivative of $\psi$ at $x$ and

(iii) $P-$ a.s. for all $x\in R$
\begin{equation}\label{ut22}
B(t,x)= \frac{1}{2}\int_0^t
\frac{(\psi_x(s,x)+\lambda(s)Y_x(s,x))'}
{Y_{xx}(s,x)}d\langle M\rangle_s(\psi_x(s,x)+\lambda(s)Y_x(s,x))
\end{equation}

{\bf Remark 3.1.} If we substitute expression of $B(t,x)$, given
by equality (\ref{ut22}), in the canonical decomposition
(\ref{ut17}) for $Y$ we obtain equation (\ref{ut20}).

{\bf Remark 3.2.}    A sufficient condition
 for twice differentiability of the value function $V(0,x)$  is given
 in Kramkov and Sirbu \cite{k-s}.

According to Proposition A1 the value process $V(t,x)$ is a
supermartingale for any $x\in R$, which admits the canonical
decomposition
\begin{equation}\label{ut23}
V(t,x)=V(0,x)+A(t,x)+\int_0^t\varphi(s,x)dM_s+m(t,x),
\end{equation}
where $-A(\cdot,x)\in{\cal A}^{+}$ and $m(\cdot,x)$ is a local
martingale strongly orthogonal to $M$ for all $x\in R_{+}$.

Assume that $V\in{\cal V}^{1,2}$. This implies that $V_x(t,x)$ is
a special semimartingale with the decomposition
\begin{equation}\label{ut24}
V_x(t,x)=V_x(0,x)+A_{(x)}(t,x)+\int_0^t\varphi_x(s,x)dM_s+m_{(x)}(t,x),
\end{equation}
where $A_{(x)}(\cdot,x)\in{\cal A}_{\text loc}$,
$m_{(x)}(\cdot,x)$ is a local martingale orthogonal to $M$ for all
$x\in R_{+}$ and $\varphi_x$ coincides with the partial derivative
of $\varphi$ ( $\mu^{K}$-a.e.). Besides
$$
A(t,x)=\int_0^t a(s,x)dK_s,
$$
for a measurable function $a(t,x)$.

Recall that the scalar product $u'\nu_tv$ for $u,v\in R^d$ we denote
by $(u,v)_{\nu_t}$.

\be{prop}\lbl{utp3.1} Assume that conditions B1), B2) are
satisfied and the value function $V(t,x)$ belongs to the class
${\cal V}^{1,2}$. Then the following inequality holds
\begin{equation}\label{ut25}
a(s,x)\le\frac{1}{2}
\frac{|\varphi_x(s,x)+\lambda(s)V_x(s-,x)|^2_{\nu_s}}{V_{xx}(s-,x)}
\end{equation}
$\text{for all}\;\;\;x\in R_+\;\;\;\;\;\mu^{K}-a.e.$ Moreover, if
the strategy $\pi^*$ is optimal then the corresponding wealth
process $X^{\pi^*}$ is a solution of the following forward SDE
\begin{equation}\label{ut26}
X_t^{\pi^*}=X_0^{\pi^*}- \int_0^t\frac{\varphi_x(s,X_s^{\pi^*})+
\lambda(s)V_x(s,X_s^{\pi^*})} {V_{xx}(s,X_s^{\pi^*})} dS_s.
\end{equation}
\ee{prop}

{\it Proof}. Using Ito-Ventzell's formula (Proposition
\ref{utp2.1}) for the function $V(t,x,\omega)\in{\cal V}^{1,2}$
and for the process $(x+\int_s^t\pi_udS_u,\;\;s\le t\le T)$ we
have
\begin{align}\label{ut27}
&V(t,x+\int_s^t\pi_udS_u)\\
\notag=& V(s,x)+\int_s^ta(u,x+\int_s^u\pi_vdS_v)dK_u\\
\notag+&\int_s^tG(u,\pi_u,x+\int_s^u\pi_vdS_v)dK_u+N_t-N_s,
\end{align}
where
$$
G(t,p,x,\omega)=V_x(t-,x)p'\nu_t\lambda(t)+p'\nu_t\vp_x(t,x)+
$$
\begin{equation}\label{ut28}
+\frac{1}{2}V_{xx}(t-,x)p'\nu_t p
\end{equation}
and $N$ is a martingale. Since by Proposition A1 of Appendix the
process \\ $(V(t,x+\int_s^t\pi_udS_u), t\in[s,T])$ is a
supermartingale for all $s\ge 0$ and $\pi\in\Pi_x$,  the process
$$
-\int_s^t\left[G(u,\pi_u,x+\int_s^u\pi_vdS_v)+
a(u,x+\int_s^u\pi_vdS_v)\right]dK_u,
$$
is increasing for any $s\ge0$.
  Hence, the process
$$
-\int_s^t\left[G(u,\pi_u,x+\int_s^u\pi_vdS_v)+
a(u,x+\int_s^u\pi_vdS_v)\right]dK^c_u,
$$
is also increasing for any $s\ge0$, where $K=K^c+K^d$ is a
decomposition of $K$ into continuous and purely discontinuous
increasing processes. Therefore, taking $\tau_s(\varepsilon)
=\inf\{t\ge s:K^c_t-K^c_s\ge\varepsilon\}$ instead of $t$ we have
that for any $\vps>0$ and $s\ge0$
\begin{align}\label{ut29}
\notag \frac{1}{\varepsilon}\int_s^{\tau_s(\varepsilon)}
a(u,x+\int_s^u\pi_vdS_v)dK^c_u\\
\le -\frac{1}{\varepsilon}\int_s^{\tau_s(\varepsilon)}
G(u,\pi_u,x+\int_s^u\pi_vdS_v)dK^c_u.
\end{align}
Passing to the limit in (\ref{ut29}) as $\vps\to 0$, from
Lemma B of \cite{MT7} we obtain that
$$
a(s,x)\le - G(s,\pi_s,x)\;\;\;\;\;\;\mu^{K^c}-a.e.
$$
for all $\pi\in\Pi$. Thus
\begin{equation}\label{utt}
a(t,x)\le \underset{\pi\in\Pi}{\essinf}
\big(-G(t,\pi_t,x)\big)\;\;\;\;\;\;\;\mu^{K^c}-a.e.
\end{equation}
On the other hand
$$
\underset{\pi\in\Pi}{\essinf}\;\big(-G(t,\pi_t,x)\big)=
\frac{|V_x(t-,x)\lambda(t)+\varphi_x(t,x)|^2_{\nu_t}}{2V_{xx}(t-,x)}
$$
$$
+\underset{\pi\in\Pi}{\essinf} \left(-\frac{1}{2}V_{xx}(t-,x)
\big|\pi_t +\frac{V_x(t-,x)\lambda(t)+
\varphi_x(t,x)}{V_{xx}(t-,x)}\big|^2_{\nu_t}\right)
$$
\begin{equation}\label{ut30}
=\frac{|V_x(t-,x)\lambda(t)+\varphi_x(t,x)|^2_
{\nu_t}}{2V_{xx}(t-,x)}.
\end{equation}
Indeed, since $V_{xx}<0$ equality  (\ref{ut30}) follows from Lemma A.1.
Thus, from ({\ref{utt}) and (\ref{ut30}) we have that for every $x\in R_+$
$$
a(t,x)\le
\frac{|V_x(t-,x)\lambda(t)+\varphi_x(t,x))|^2_{\nu_t}}{2V_{xx}(t-,x)},
\;\;\;\;\mu^{K^c}\;\;\;a.e.
$$
Since $\mu^K$- a.e. $a(t,x)\ge0$ and $\mu^{K^d}\{\nu\ne 0\}=0$ we
obtain that
\begin{equation}\label{ut31}
a(t,x)\le
\frac{|V_x(t-,x)\lambda(t)+\varphi_x(t,x)|^2_{\nu_t}}{2V_{xx}(t-,x)},
\;\;\;\;\mu^{K}\;\;\;a.e.
\end{equation}

Conditions C2) and C5) imply that inequality (\ref{ut31}) holds
$\mu^K$-a.e. for all $x\in R$.

Let us show now that if the  strategy $\pi^*$ is optimal then the
corresponding wealth process $X^{\pi^*}$ is a solution of equation
(\ref{ut26}). Let $\pi^*(s,x)$ be the optimal strategy and denote
by $X^*_t(s,x)=x+\int_s^t\pi_u^*(s,x)dS_u$ the corresponding
wealth process.

By the optimality principle the process
$V(t,x+\int_s^t\pi_u^*(s,x)dS_u)$ is a martingale on the time
interval $[s,T]$ and the Ito-Ventzell formula implies that
$\mu^{K}$-a.s.
$$
a(t,X_t^*(s,x))+(\lambda_t,\pi_t(s,x))_{\nu_t}V_x(t-, X_t^*(s,x))+
$$
\begin{equation}\label{ut32}
(\vp_x(t,X_t^*(s,x)),\pi^*_t(s,x))_{\nu_t}
+\frac{1}{2}|\pi^*_t(s,x)|^2_{\nu_t}V_{xx}(t-,X_t^*(s,x))=0.
\end{equation}
It follows from (\ref{ut31}) and  (\ref{ut32}) that $\mu^{K}$-a.e.
$$
V_{xx}(t-,X_t^*(s,x))
\big|\pi_t^*(s,x)+\frac{\varphi_x(t,X^*_t(s,x))+
\lambda(t)V_x(t-,X_t^*(s,x))}
{V_{xx}(t-,X_t^{*}(s,x))}\big|^2_{\nu_t} \ge 0.
$$
Since  $V_{xx}<0$, integrating the latter relation by $dK_u$ we
obtain that
$$
\int_s^t \big(\pi_u^*(s,x)+\frac{\varphi_x(u,X^*_u(s,x))+
\lambda(u)V_x(u,X_u^*(s,x))}
{V_{xx}(u,X_u^{*}(s,x))}\big)'d\langle M\rangle_u\times
$$
\begin{equation}\label{ut33}
\times\big(\pi_u^*(s,x)+\frac{\varphi_x(u,X^*_u(s,x))+\lambda(u)V_x(u,X_u^*(s,x))}
{V_{xx}(u,X_u^{*}(s,x))}\big)=0.
\end{equation}
The Kunita--Watanabe inequality and (\ref{ut33}) imply that the
semimartingale
$$
\int_s^t\big(\pi^*_u(s,x)+\frac{\varphi_x(u,X^*_u(s,x))+
\lambda(u)V_x(u,X^*_u(s,x))} {V_{xx}(u,X^*_u(s,x))}\big) dS_u
$$
is indistinguishable from zero (since its ${\cal S}^2$-norm is
zero) and we obtain that the wealth process of $\pi^*$ satisfies
equation
\begin{equation}\label{ut34}
X^*_t(s,x)=x- \int_s^t\frac{\varphi_x(u,X^*_u(s,x))+
\lambda(u)V_x(u,X^*_u(s,x))} {V_{xx}(u,X^*_u(s,x))} dS_u
\end{equation}
which gives equation (\ref{ut26}) for $s=0$.\qed


Recall that the process $Z$ belongs to the class $D$ if the family
of random variables $Z_\tau I_{(\tau\le T)}$ for all stopping
times $\tau$ is uniformly integrable.

Under additional condition

\noindent {\bf C*)} $(X_t^*(s,x),\;t\ge s)$ is a continuous
function of $(s,x)$ $P-$a.s. for each $t\in [s,T]$,

we shall show
that the value function $V$ satisfies equation
(\ref{ut20})-(\ref{ut21}).

This condition is satisfied, e.g., if  the optimal wealth
process\\
$(X_t^*(s,x),\;t\ge s)$ does not depend on $s$ and $x$, which we
have in cases of power, logarithmic and exponential utility
functions.

\be{th}\lbl{utt3.1} Let $V\in{\cal V}^{1,2}$ and assume that
conditions B1)-B3), C*) are satisfied. Then the value function is
a solution of BSPDE (\ref{ut20})-(\ref{ut21}), i.e.,
$$
V(t,x)= V(0,x)
+\frac{1}{2}\int_0^t \frac{(\varphi_x(s,x)+\lambda(s)V_x(s,x))'}
{V_{xx}(s,x)}d\langle M\rangle_s(\varphi_x(s,x)+\lambda(s)V_x(s,x))
$$
\begin{equation}\label{ut35}
+\int_0^t\varphi(s,x)dM_s+m(t,x),\;\;\; V(T,x)=U(x).
\end{equation}
Moreover, the strategy $\pi^*$ is optimal if and only if the
corresponding wealth process $X^{\pi^*}$ is a solution of the
forward SDE (\ref{ut26}), such that the process $V(t,X^{\pi^*})$
is from the class $D$. \ee{th}

{\it Proof.} Let $\pi^*(s,x)$ be the optimal strategy. By
optimality principle $(V(t,X_t^*(s,x)),t\ge s)$ is a martingale.
Therefore, using Ito-Ventzell's formula, taking (\ref{ut33}) in
mind, we have
$$
\int_s^t\big[a(u,X_u^*(s,x))- g(u,X_u^*(s,x))+
$$
$$
+\big|\pi_u^*(s,x) + \frac{V_x(u,X_u^*(s,x))\lambda(u)+
\varphi_x(u,X_u^*(s,x))}{V_{xx}(u,X_u^*(s,x))}\big|^2_{\nu_u}
\big]dK_u=0,
$$
$$
\text{for all}\;\;t\ge s \;\;P-a.s.,
$$
where
$$
g(s,x)=\frac{1}{2}
\frac{|\varphi_x(s,x)+\lambda(s)V_x(s,x)|^2_{\nu_s}}
{V_{xx}(s,x)}.
$$
It follows from (\ref{ut33}) that $\mu^K-a.e.$
$$
\big|\pi_u^*(s,x) + \frac{V_x(u,X_u^*(s,x))\lambda(u)+
\varphi_x(u,X_u^*(s,x))}{V_{xx}(u,X_u^*(s,x))} \big|^2_{\nu_u}=0
$$
and by (\ref{ut25})
\begin{equation}\label{ut36}
a(s,x)\le g(s,x)\;\;\;\;\mu^K-a.e.
\end{equation}
Thus,
$$
\int_s^t[a(u,X_u^*(s,x))- g(u,X_u^*(s,x))]dK_u=0,\; t\ge s
\;\;\;\;\; P- a.s.
$$
This implies that $(a(s,x)-g(s,x))(K_s-K_{s-})=0$ for any
$s\in[0,T]$. Therefore,
\begin{equation}\label{ut37}
a(s,x)=g(s,x)\;\;\;\;\mu^{K^d}-a.e.
\end{equation}

On the other hand
$$
\int_0^T\frac{1}{\vps}\int_s^{\tau^\vps_s}[a(u,X_u^*(s,x))-
g(u,X_u^*(s,x))]dK^c_udK^c_s=0,\; \;\;P-\;a.s.
$$
and by Proposition B of \cite{MT7} we obtain that
$$
\int_0^T[a(s,x)-g(s,x)]dK^c_s=0,\;\;P-a.s.
$$
Now (\ref{ut36}), (\ref{ut37}) and the latter relation result
equality $a(s,x)=g(s,x)\;\;\mu^K-a.e.$, hence
$$
A(t,x)=\frac{1}{2}\int_0^t
\frac{(\varphi_x(s,x)+\lambda(s)V_x(s,x))'}
{V_{xx}(s,x)}d\langle M\rangle_s(\varphi_x(s,x)+\lambda(s)V_x(s,x))
$$
and $V(t,x)$ satisfies (\ref{ut20})-(\ref{ut21}).

If $\hat\pi$ is a strategy such that the corresponding wealth
process $X^{\hat\pi}$ satisfies equation (\ref{ut26}) and
$V(t,X_t^{\hat\pi})$ is from the class $D$, then $\hat\pi$ is
optimal. Indeed, using the Ito-Ventzell formula and equations
(\ref{ut26}) and (\ref{ut35}) we obtain that $V(t,X_t^{\hat\pi})$
is a local martingale, hence a martingale, since it belongs to the
class $D$. Therefore $\hat\pi$ is optimal by optimality
principle.\qed

{\bf Remark 3.3.} For the utility functions which take finite
values on all real line  Proposition 3.1 and Theorem 3.1 are also
true if we  choose a suitable class of trading strategies. E.g.,
let $\Pi_x$ be one of the class introduced by  Schachermayer
(2003)(${\cal H}_i(x)$ or ${\cal H}^{'}_i(x)$, for $i=1, 2$ or 3).
The proof of abovementioned assertions is the same, there is a
minor difference  only in the proof of equality (\ref{ut30}),
where instead of Lemma A.1 the following arguments  should be
used: Conditions (C3), (C4) and (C6) imply that for each $x$ the
process
$$
\frac{V_x(t-,x)\lambda(t)+
\varphi_x(t,x)}{V_{xx}(t-,x)}
$$
is predictable and $S$-integrable. Therefore, there exists a sequence of stopping times
$(\tau_n(x), n\ge1)$ with $\tau_n(x)\uparrow T$ for all $x\in R$ such that the wealth process
corresponding to the strategy
$$
\pi^n_t=-I_{[0,\tau_n]}\frac{V_x(t-,x)\lambda(t)+
\varphi_x(t,x)}{V_{xx}(t-,x)}
$$
is bounded and hence $\pi^n\in\Pi$ for each $n$. Therefore,
$$
0\le\underset{\pi\in\Pi}{\essinf} \left(-\frac{1}{2}V_{xx}(t-,x)
\big|\pi_t +\frac{V_x(t-,x)\lambda(t)+
\varphi_x(t,x)}{V_{xx}(t-,x)}\big|^2_{\nu_t}\right)
$$
$$
\le\left(-\frac{1}{2}V_{xx}(t-,x)
\big|\frac{V_x(t-,x)\lambda(t)+
\varphi_x(t,x)}{V_{xx}(t-,x)}\big|^2_{\nu_t}\right)I_{(\tau_n(x)\le t)}\to 0,\;\;\;
\mu^{K^c}\;\;\;\text{a.s.},
$$
which implies equality (\ref{ut30}).

\be{th}\lbl{utt3.2} Let conditions B1)-B3) be satisfied. If the
pair $(Y,\mathcal X)$ is a solution of the Forward-Backward
Equation
\begin{align}\label{ut38}
Y&(t,x)=U(x)\\
\notag &-\frac{1}{2}\int_t^T \frac{((\psi_x(s,x)+\lambda(s)Y_x(s,x))'}
{Y_{xx}(s,x)}d\langle M\rangle_s(\psi_x(s,x)+\lambda(s)V_x(s,x))
\\
\notag &-\int_t^T\psi(s,x)dM_s+L(T,x)-L(t,x)
\end{align}
\begin{equation}\label{ut39}
{\mathcal X}_t=x- \int_0^t\frac{\psi'_x(s,{\mathcal X}_s)+
Y_x(s,{\mathcal X}_s)\lambda(s)} {Y_{xx}(s,{\mathcal X}_s)} dS_s,
\end{equation}
${\mathcal X}\ge 0,\;\;Y\in{\cal V}^{1,2}$ and  $Y(t,{\cal X}_t)$ belongs to
the class $D$, then such solution is unique.
\ee{th}

{\it Proof}. Using the Ito-Ventzell's formula for
$Y(t,x+\int_s^t\pi_udS_u)$ we have
\begin{align}\label{ut40}
&Y(t,x+\int_s^t\pi_udS_u)\\
\notag
&=Y(s,x)+\int_ s^tb(u,x+\int_s^u\pi_vdS_v)dK_u\\
\notag
&+\int_s^tG(u,\pi_u,c+\int_s^u\pi_vdS_v)dK_u+N_t-N_s,
\end{align}
where
$$
G(t,p,x,\omega)=Y_x(t-,x)p'\nu_t\lambda(t)+p'\nu_t\psi_x(t,x)+
\frac{1}{2}Y_{xx}(t-,x)p'\nu_t p
$$
and $N$ is a local martingale.

Since $Y$ solves (\ref{ut38}), then equality (\ref{ut22}) is
valid, which implies that $Y(t,x+\int_s^t\pi_udS_u)$ is a local
supermartingale for each $\pi\in\Pi$.

Let $\tau_n=\inf\{t:Y(t,x+\int_s^t\pi_udS_u)\ge n\}\wedge T$.
\footnote{$)$\,It is assumed that $\inf\emptyset=\infty$ and
$a\wedge b$ denotes $\min\{a,b\}$ }$)$
Then by supermartingale property and the monotone convergence theorem we have
$$
Y(s,x)\ge E\big(Y(\tau_n,x+\int_s^{\tau_n}\pi_udS_u)|{\cal F}_s\big)
$$
$$
\ge E\big(n\wedge U(x+\int_s^T\pi_udS_u)|{\cal F}_s\big)
\overset{n\to\infty}\longrightarrow
E\big(U(x+\int_s^T\pi_udS_u)|{\cal F}_s\big).
$$
i.e.
$$
Y(s,x)\ge E\big(U(x+\int_s^T\pi_udS_u)|{\cal F}_s\big),\;\;\forall\pi\in\Pi_x,
$$
which implies  that
\begin{equation}\label{ut41}
Y(s,x)\ge V(s,x).
\end{equation}
Using now the Ito-Ventzell's formula for $Y(t,{\mathcal X}_t)$
taking into account that $Y$ satisfies (\ref{ut38}) and $\mathcal
X$ solves (\ref{ut39}) we obtain that $Y(t,{\mathcal X}_t)$ is a
local martingale and, hence, it is a martingale, since
$Y(t,{\mathcal X}_t)$ is from the class $D$.

Therefore, since ${\mathcal X}_0=x,\;Y(T,x)=U(x)$ we have that

\begin{equation}\label{ut42}
Y(t,x)=E\big(U(x-\int_t^T\frac{Y_x(u,{\mathcal X}_u)\lambda_u+
\psi_x(u,{\mathcal X}_u)}{Y_{xx}(u,{\mathcal X}_u)}dS_u)/{\cal
F}_t\big).
\end{equation}
 Since $-\frac{\lambda(u)Y_x(u,{\mathcal X}_u)+
\psi_x(u,{\mathcal X}_u)}{Y_{xx}(u,{\mathcal X}_u))}\in\Pi_x$,
from (\ref{ut41}) and (\ref{ut42}) we obtain that
\begin{equation}\label{ut43}
Y(t,x)=V(t,x),
\end{equation}
hence solution of (\ref{ut38}) is unique if it exists and coincides with
the value function. This implies that under conditions of theorem
$V\in {\cal V}^{1,2}$.

Therefore, it follows from (\ref{ut43}) and (\ref{ut39})
that $\mathcal X$ satisfies  equation (\ref{ut26}). Besides,
according to Proposition \ref{utp3.1} the solution of (\ref{ut26})
is the optimal wealth process, hence $\mathcal X =X^{\pi*}$ by the
uniqueness of the optimal strategy for the problem (\ref{ut3})
(see Remark 2.2).

\

\section{Utility maximization problem for
power, logarithmic and exponential utility functions}

\

In this section we calculate the value function and give
constructions of optimal strategies for the utility maximization
problem corresponding to the cases of power, logarithmic and
exponential utility functions.

{\it Power Utility.}

Let $U(x)=x^p,\; p\in (0,1)$. Then (\ref{ut3}) corresponds to power
utility maximization problem
\begin{equation}\label{pw}
\text{to maximize}\;\;\;\;\;
E(x+\int_0^T\pi_udS_u)^p\;\;\;\;\text{over all}\;\;\;\;\pi\in\Pi_x
\end{equation}
where $\Pi_x$ is a class of admissible strategies.

In this case the value function $V(t,x)$ is of the form $x^pV_t$,
where $V_t$ is a special semimartingale.
Indeed, since $\Pi_x$ is a cone (for any $x>0$ the strategy $\pi$
belongs to $\Pi_x$ iff $\frac{\pi}{x}\in\Pi_1$) , we have
$$
V(t,x)=\underset{\pi\in\Pi_x}\esssup E\big((x+\int_t^T\pi_udS_u
)^p/{\mathcal F}_t\big)
$$
$$
=x^p\underset{\pi\in\Pi_x}\esssup E\big((1+\int_t^T\frac{\pi_u}{x}dS_u
)^p/{\mathcal F}_t\big)=x^pV_t,
$$
where
$$
V_t=\underset{\pi\in\Pi_1}\esssup E\big((1+\int_t^T\pi_udS_u
)^p/{\mathcal F}_t\big)
$$
is a supermartingale by optimality principle.

Let $V_t=V_0+A_t+N_t$ be the canonical decomposition of $V_t$,
where $A$ is a decreasing process and $N$ is a local martingale.
Using the G--K--W decomposition we have that
\begin{equation}
 V_t=V_0+A_t+\int_0^t\varphi_sdM_s + L_t,
\end{equation}
where $L$ is a local martingale with $<L,M>=0$.

It is evident that for $U(x)=x^p$ the condition  (\ref{ae}) is
satisfied and the optimal strategy for the problem (\ref{pw})
exists.  Since in this case $V(t,x)=x^pV_t$ it is also evident that
$V(t,x)\in{\mathcal V}^{1,2}$ and all conditions of Theorem 3.1 are satisfied
(note that one can take $-A+\langle M\rangle$ as a dominated process $K$).

Therefore we have the following corollary of Theorem \ref{utt3.1}

\be{th}\lbl{utt4.1} If  $U(x)=x^p, p\in (0,1)$, then the value
function $V(t,x)$ is of the form $x^pV_t$, where $V_t$ satisfies
the following backward stochastic differential equation (BSDE)
$$
V_t= V_0 +\frac{q}{2}\int_0^t \frac{(\varphi_s+\lambda_sV_s)'}
{V_s}d\langle M\rangle_s(\varphi_s+\lambda_sV_s)
$$
\begin{equation}\label{power}
+\int_0^t\varphi_sdM_s+L_t,\;\;V_T=1,
\end{equation}
where $q=\frac{p}{p-1}$ and $L$ is a local martingale strongly
orthogonal to $M$.

Besides, the optimal wealth process is a solution of the linear
equation
\begin{equation}
X_t^{*}=x-(q-1)\int_0^t\frac{\varphi_u+\lambda_uV_u}
{V_u}X^{*}_udS_u
\end{equation}
Therefore,
$$
X_t^{*}=x{\mathcal E}_t(-(q-1)(\frac{\varphi}{V}+\lambda)\cdot S)
$$
and the optimal strategy is of the form
$$
\pi_t^*=-x(q-1){\mathcal
E}_t(-(q-1)(\frac{\varphi}{V}+\lambda)\cdot S)
(\frac{\varphi_t}{V_t}+\lambda_t).
$$
\ee{th}

{\bf Remark 4.1.} If there is a martingale measure $Q$ that satisfies the
Muckenhoupt condition $A_\alpha(P)$ for $\alpha=\frac{1}{p}$ then the process $V$ is bounded.
Indeed,  by the H\"older inequality for any $\pi\in\Pi$
$$
E\big((1+\int_t^T\pi_udS_u )^p/{\mathcal F}_t\big)=
E^Q\big((1+\int_t^T\pi_udS_u )^p\frac{Z_t^Q}{Z_T^Q}|{\mathcal
F}_t\big)\le
$$
$$
\le \big(E^Q((1+\int_t^T\pi_udS_u)|{\mathcal F}_t)\big)^p
 \big(E^Q\big(\big({Z_t^Q}/{Z_T^Q}\big)^{\frac{1}{1-p}}|{\mathcal F}_t\big)\big)^{1-p}\le
$$
$$
\le
\big(E\big(\big({Z_t^Q}/{Z_T^Q}\big)^{\frac{1}{\frac{1}{p}-1}})|
{\mathcal F}_t\big)\big)^{1-p}\le C^{1-p}.
$$
Under condition $A_{\frac{1}{p}}(P)$ equation
(\ref{power}) admits a unique bounded strictly positive solution. This follows from Theorem 3.2,
since in this case the process ${\mathcal E}_t(-(q-1)(\frac{\varphi}{V}+\lambda)\cdot S)$
belongs to the class $D$.

Now we shall consider two cases when equation (\ref{power}) admits
an explicit solution

Case 1.

Let $S_t(q)= M_t+q\int_0^td\langle M\rangle_s\lambda_s$ and let $Q(q)$ be a
measure defined by $dQ(q)={\mathcal E}_T(-q\lambda\cdot M)dP$.
Note that $S(q)$ is a local martingale under $Q(q)$ by Girsanov's
theorem.

Assume that
\begin{equation}\label{4.5}
e^{\frac{q(q-1)}{2}\langle\lambda\cdot M\rangle_T}=c+\int_0^Th_udS_u(q),
\end{equation}
where $c$ is a constant and $h$ is a predictable $S(q)$-integrable
process such that $h\cdot S(q)$ is a $Q(q)$-martingale.

This condition is satisfied iff the $q$-optimal martingale measure
coincides with the minimal martingale measure. For diffusion
market models this condition is fulfilled for so called "almost
complete" models, i.e., when the market price of risk is
measurable with respect to the filtration generated by price
processes of basic securities.

Let condition (\ref{4.5}) be satisfied. Let us consider the
process
\begin{equation}\label{4.6}
Y_t=\big(E({\mathcal E}^q_{t,T}(-\lambda\cdot
M)/F_t)\big)^{\frac{1}{1-q}}.
\end{equation}
Since
$$
{\mathcal E}_t^q(-\lambda\cdot M)= {\mathcal E}_t(-q\lambda\cdot
M) e^{\frac{q(q-1)}{2}\langle\lambda\cdot M\rangle_t},
$$
condition (\ref{4.5}) implies that
$$
Y_t=\big(E^{Q(q)}(e^{\frac{q(q-1)}{2} (\langle\lambda\cdot M\rangle_T-\langle\lambda\cdot M\rangle_t}/F_t)\big)^{\frac{1}{1-q}}=
$$
$$
= e^{\frac{q}{2}\langle\lambda\cdot M\rangle_t}
\big(c+\int_0^th_udS_u(q)\big)^{\frac{1}{1-q}}.
$$
By the It\^o formula
$$
Y_t=Y_0+\frac{q}{2}\int_0^tY_s\lambda'_sd\langle M\rangle_s\lambda_s+
\frac{q}{1-q}\int_0^t\frac{Y_s\lambda'_s}{c+(h\cdot S(q))_s}
d\langle M\rangle_sh_s
$$
\begin{equation}\label{4.7}
+\frac{q}{2}\frac{1}{(1-q)^2}\int_0^t\frac{Y_sh'_s}{(c+(h\cdot
S(q))_s)^2} d\langle M\rangle_sh_s+
\frac{1}{1-q}\int_0^t\frac{Y_sh_s}{c+(h\cdot S(q))_s} dM_s
\end{equation}
and denoting   $\frac{1}{q-1}\frac{Y_sh_s}{c+(h\cdot S(q))_s}$ by
$\psi_s$ we obtain that
$$
Y_t= Y_0+\frac{q}{2}\int_0^t \frac{(\psi_s+\lambda_sY_s)'}
{Y_s}d\langle M\rangle_s(\psi_s+\lambda_sY_s)+\int_0^t\psi_sdM_s.
$$
It is evident from (\ref{4.6}) that $Y_T=1$. Thus the triple
$(Y,\psi, L)$, where $\psi=\frac{1}{q-1}\frac{Yh}{c+h\cdot S(q)}$,
$L=0$ and $Y$ defined by (\ref{4.6}), satisfies equation
(\ref{power}).

Case 2.

Assume that
\begin{equation}\label{4.88}
e^{-\frac{q}{2}\langle\lambda\cdot M\rangle_T}=c+m_T,
\end{equation}
where $c$ is a constant and $m$ is a martingale strongly
orthogonal to $M$.

For diffusion market models this condition is satisfied when the
market price of risk is measurable with respect to the filtration
independent relative to the asset price process.

Let us consider the process
\begin{equation}\label{4.9}
Y_t=E(e^{-\frac{q}{2}(\langle\lambda\cdot M\rangle_T-\langle\lambda\cdot M\rangle_t)}/F_t).
\end{equation}

Condition (\ref{4.88}) implies that
$$
Y_t=e^{\frac{q}{2}\langle\lambda\cdot M\rangle_t}(c+m_t)
$$
and by the It\^o formula
$$
Y_t=Y_0+\frac{q}{2}\int_0^tY_sd\langle\lambda\cdot M\rangle_s+ \int_0^t
e^{\frac{q}{2}\langle\lambda\cdot M\rangle_s}dm_s.
$$
It follows from here that the triple $(Y,\psi, L)$, where
$\psi=0$, $L_t=\int_0^t e^{\frac{q}{2}\langle\lambda\cdot M\rangle_s}dm_s$ and
$Y$ defined by (\ref{4.9}), satisfies equation (\ref{power}) and the optimal strategy
is
$$
\pi_t^*=x(1-q)\lambda_t{\mathcal
E}_t((1-q)\lambda\cdot S).
$$

 {\it Logarithmic Utility}

For the logarithmic utility
$$
U(x)= \log x,\;\;\;x>0
$$
the value function of corresponding utility maximization problem
takes the form
$$
V(t,x)=\log x + V_t,
$$
where $V_t$ is a special semimartingale.

Indeed, since for any $x>0$ the strategy $\pi$ belongs to $\Pi_x$
iff $\frac{\pi}{x}\in\Pi_1$, we have
$$
V(t,x)=\underset{\pi\in\Pi_x}{\text{\rm esssup}}\text{ }
E\big(\log\big(x+\int_t^T\pi_udS_u \big)/{\mathcal F}_t\big)
$$
$$
=\underset{\pi\in\Pi_x}{\text{\rm esssup}}\text{ } E\big(\log
x\big(1+\int_t^T\frac{\pi_u}{x}dS_u\big) /{\mathcal F}_t\big)
$$
$$
=\log x +\underset{\pi\in\Pi_x}{\text{\rm esssup}}\text{ } E\big(\log
\big(1+\int_t^T\frac{\pi_u}{x}dS_u \big)/{\mathcal F}_t\big)= \log x+V_t,
$$
where
$$
V_t=\underset{\pi\in\Pi_1}{\text{\rm esssup}}\text{ }
E\big(\log\big(1+\int_t^T\pi_udS_u \big)/{\mathcal F}_t\big)
$$
is a supermartingale by the optimality principle.

It is also evident that  all conditions of Theorem \ref{utt3.1}
are satisfied. In this case
$\varphi_x(t,x)=0,V_x(t,x)=\frac{1}{x},V_{xx}(t,x)=-\frac{1}{x^2}$
and equation (\ref{ut26}) gives the following expression for
$V_t$
$$
V_t= V_0-\frac{1}{2}\langle\lambda\cdot M\rangle_t+\int_0^t\varphi_sdM_s+L_t,
$$
\begin{equation}
V_T=0,
\end{equation}
which admits an explicit solution
$$
V_t=-\frac{1}{2}E(\langle\lambda\cdot M\rangle_T-\langle\lambda\cdot M\rangle_t/F_t).
$$
Thus, we have the following corollary of Theorem \ref{utt3.1}
\be{th}\lbl{utt4.3} If  $U(x)=\log x$, then the value function of
the problem is represented as
$$
V(t,x)= \log x - \frac{1}{2}E(\langle\lambda\cdot M\rangle_T-\langle\lambda\cdot M\rangle_t/F_t).
$$
Besides, the optimal wealth process is a solution of the linear
equation
\begin{equation}
X_t^{*}=x+\int_0^t\lambda_u X^*_udS_u.
\end{equation}
Thus,
$$
X_t^{*}=x{\mathcal E}_t(\lambda\cdot S)
$$
and the optimal strategy is of the form
$$
\pi_t^*=\lambda_tX_t^*=x\lambda_t{\mathcal E}_t(\lambda\cdot S).
$$
\ee{th}

{\it Exponential Utility}

Let us consider the case of exponential utility function
$$
U(x)=-e^{-\gamma (x-H)}
$$
with risk aversion parameter $\gamma\in(0,\infty)$, where $H$ is a
bounded contingent claim describing a random payoff at time $T$. We assume
that $H$ is bounded $F_T$-measurable random variable.

For any $Q\in{\cal M}^e$ let $(Z_t^Q,t\in[0,T])$ be the density process of $Q$
with respect to $P$ and assume that
$$
{\cal M}^e_{\ln}=\{Q\in{\cal M}^e:EZ^Q_T\ln Z^Q_T<\infty\}\neq\emptyset.
$$
We define the space of trading strategies $\Pi$  as the space of all predictable $S$-integrable
processes $\pi$ such that the corresponding wealth process $X^\pi$ is a martingale relative to
any $Q\in{\cal M}^e_{\ln}$. So, $\Pi$ is the space $\Theta_2$ from Delbaen et al. (2002) and the space
${\cal H}_2$ from Schachermayer (2003).

Let us consider the maximization problem
\begin{equation}\label{exp}
\max_{\pi\in\Pi} E(-e^{-\gamma(x+\int_0^T\pi_udS_u-H)}),
\end{equation}
the maximal expected utility we can achieve by starting with
initial capital $x$, using some strategy $\pi\in\Pi$ and paying
out $H$ at time $T$.

The corresponding  value function
\begin{equation}\label{vtx}
V(t,x)=\underset{\pi\in\Pi_x}{\text{\rm esssup}}\text{ }
E(-e^{-\gamma(x+\int_t^T\pi_udS_u-H)}/{\mathcal F}_t)
\end{equation}
is of the form $V(t,x)=-e^{-\gamma x}V_t$, where

\begin{equation}\label{utve}
V_t=\underset{\pi\in\Pi_x} {\text{\rm
essinf}}\text{ } E(e^{-\gamma(\int_t^T\pi_udS_u-H)}|{\mathcal
F}_t) \end{equation}

is a special semimartingale.

Let $V_t=V_0+A_t+N_t$ be the canonical decomposition of $V_t$,
where $A$ is an increasing process and $N$ is a local matingale.
Using the G--K--W decomposition we have that
\begin{equation}
 V_t=V_0+A_t+\int_0^t\varphi_sdM_s + L_t,
\end{equation}
where $L$ is a local martingale with $\la L,M\ra=0$.

Since ${\cal M}^e_{\ln}\neq\emptyset$, the optimal strategy in the class $\Pi$ exists
and $V_t >0$ for all $t$ (see, e.g., Delbaen at al. (2002) and
Yu. Kabanov and Ch. Stricker (2002)). It is evident that
$V(t,x)=-e^{-\gamma x}V_t\in{\mathcal V}^{1,2}$  and
all conditions of Theorem 3.1 are satisfied.

Therefore,  Theorem \ref{utt3.1} and Remark 3.3 imply the validity of the
following

\be{th}\lbl{utt4.2} The value function (\ref{vtx}) is of the form
$-e^{-\gamma x}V_t$, where $V_t$ satisfies the BSDE
\begin{equation}\lbl{bexp}
V_t= V_0+\frac{1}{2}\int_0^t \frac{(\varphi_s+\lambda_sV_s)^2}
{V_s}d\langle M\rangle_s
+\int_0^t\varphi_sdM_s+L_t
\end{equation}
with the boundary condition
$$
V_T=e^{\gamma H}.
$$
where  $L$ is a local martingale strongly orthogonal to $M$.

Besides, the optimal wealth process is expressed as
\begin{equation}\label{frexp}
X_t^{*}=x+\int_0^t\frac{\varphi_u+\lambda_uV_u} {\gamma V_u}dS_u
\end{equation}
 and the optimal strategy is of the form
$$
\pi_t^*=\frac{\varphi_t+\lambda_t V_t}{\gamma V_t}.
$$
\ee{th}

{\bf Remark 4.2} It is evident that $V_t\le E(e^{\gamma H}|F_t)\le const$.
If there exists a martingale measure $Q$ that satisfies the Reverse H\"older
$R_{Llog L}$ condition, i.e., if
$$
E(\frac{Z^Q_T}{Z^Q_t}\ln\frac{Z^Q_T}{Z^Q_t}|F_t)\le C
$$
for all $t$,  then there is a constant $c>0$ such that $V_t\ge c$
for all $t$ (see \cite{D-Gr-Rh-S-S-S}, \cite{M-sch}). Under   $R_{Llog L}$
condition  the value process $V$ is the unique bounded strictly positive solution
of BSDE (\ref{bexp}).

Now we shall give explicit solutions of equation
(\ref{bexp}) in two extreme cases.

Case 1.

Assume that
\begin{equation}\label{4.5ex}
 \gamma H-\frac{1}{2}\langle\lambda\cdot M\rangle_T=c+\int_0^Th_udS_u,
\end{equation}
where $c$ is a constant and $h$ is a predictable $S$-integrable
process such that $h\cdot S$ is a martingale with respect to the
minimal martingale measure.

This condition is satisfied iff the minimal entropy martingale  measure
coincides with the minimal martingale measure and $H$ is
attainable. For diffusion market models this condition is
fulfilled for so called "almost complete" models, i.e., when the
market price of risk is measurable with respect to the filtration
generated by price processes of basic securities.

 Similarly to the case of power utility one can show that
 the triple $(Y,\psi, L)$, where
 $$
 Y_t=e^{E^{Q^{min}}( \gamma H-\frac{1}{2}<\lambda\cdot
 M>_{tT}/F_t)},\;\;\;\;\psi_t=Y_th_t,\;\;\;L_t=0
 $$
 satisfies equation (\ref{bexp}) and the optimal strategy is
 $$
 \pi^*_t=\frac{1}{\gamma}(\lambda_t+h_t).
 $$

Case 2.

Assume that
\begin{equation}\label{4.8}
e^{\gamma H-\frac{1}{2}\langle\lambda\cdot M\rangle_T}=c+m_T,
\end{equation}
where $c$ is a constant and $m$ is a martingale strongly
orthogonal to $M$.

For diffusion market models this condition is satisfied when the
market price of risk and $H$ are measurable with respect to the filtration
independent relative to the asset price process.

One can show that the triple $(Y,\psi, L)$, where
$$
Y_t=e^{E( \gamma H-\frac{1}{2}\langle\lambda\cdot M\rangle_{tT}/F_t)},\;\;\;\;\psi_t=0,\;\;\;L_t=\int_0^te^{\frac{1}{2}\langle\lambda\cdot M\rangle_s}dm_s
 $$
 satisfies equation (\ref{bexp}) and the optimal strategy is $\pi^*_t=\frac{1}{\gamma}\lambda_t$.

\

{\it Quadratic Utility.}

Let $U(x)=2bx-x^2$, where $b$ is a positive constant.

Assume that
$$
{\cal M}^e_{2}=\{Q\in{\cal M}^e:E(Z^Q_T)^2<\infty\}\neq\emptyset.
$$
and let $\Pi$  be the space of all predictable $S$-integrable
processes $\pi$ such that $\int_0^T\pi_udS_u$ is in $L^2(P)$ and
the stochastic integral $\int_0^t\pi_udS_u$ is a martingale relative to
any $Q\in{\cal M}^e_2$.

In this case (\ref{ut3}) corresponds to the
utility maximization problem
\begin{equation}\label{qua}
\text{to maximize}\;\;\;\;\;
E[x+2b\int_0^T\pi_udS_u-(\int_0^T\pi_udS_u)^2]\;\;\;\;\text{over all}\;\;\;\;\pi\in\Pi,
\end{equation}
which is equivalent to the problem
\begin{equation}\label{qmin}
\text{to minimize}\;\;\;\;\;
E(x+\int_0^T\pi_udS_u-b)^2\;\;\;\;\text{over all}\;\;\;\;\pi\in\Pi.
\end{equation}
This is the mean variance
hedging problem with a constant contingent claim.

In this case the value function of (\ref{qua}) is of the form $V(t,x)=b^2-(x-b)^2V_t$,
where
$$
V_t=\underset{\pi\in\Pi}\essinf E\big((1+\int_t^T\pi_udS_u
)^2/{\mathcal F}_t\big)
$$
is a supermartingale by optimality principle.

Let $V_t=V_0+A_t+N_t$ be the canonical decomposition of $V_t$,
where $A$ is an increasing process and $N$ is a local martingale.
Using the G--K--W decomposition we have that
\begin{equation}
 V_t=V_0+A_t+\int_0^t\varphi_sdM_s + L_t,
\end{equation}
where $L$ is a local martingale with $\langle L,M\rangle=0$.

Since ${\cal M}^e_{2}\neq\emptyset$, the optimal strategy in the class $\Pi$ exists
and $V_t >0$ for all $t$ (see, e.g., Gourieroux et al. 1998 or Heath et al. 2001) . It is evident that
$V(t,x)=b^2-(x-b)^2V_t$ belongs to the class ${\mathcal V}^{1,2}$  and
all conditions of Theorem 3.1 are satisfied
(again one can take $A+\langle M\rangle$ as a dominated process $K$)

Therefore,  Theorem \ref{utt3.1} and Remark 3.3  imply the followin assertion

\be{th}\lbl{utt4.1} If  $U(x)=2bx-x^2, b\ge0$, then the value
function $V(t,x)$ is of the form $b^2-(x-b)^2V_t$, where $V_t$ satisfies
the BSDE
$$
V_t= V_0 +\int_0^t \frac{(\varphi_s+\lambda_sV_s)'}
{V_s}d\langle M\rangle_s(\varphi_s+\lambda_sV_s)
$$
\begin{equation}\label{power2}
+\int_0^t\varphi_sdM_s+L_t,\;\;V_T=1.
\end{equation}

Besides, the optimal wealth process is a solution of the linear
equation
\begin{equation}
X_t^{*}=x-\int_0^t\frac{\varphi_u+\lambda_uV_u}
{V_u}X^{*}_udS_u
\end{equation}
Therefore,
$$
X_t^{*}=x{\mathcal E}_t(-(\frac{\varphi}{V}+\lambda)\cdot S)
$$
and the optimal strategy is of the form
$$
\pi_t^*=-x{\mathcal
E}_t(-(\frac{\varphi}{V}+\lambda)\cdot S)
(\frac{\varphi_t}{V_t}+\lambda_t).
$$
\ee{th}

\

\section{Diffusion market models}

\

The main task of this section is to establish a connection between
the semimartingale backward equation for the value process and the
classical Bellman equation for the value function related to the
utility maximization problem in the case of Markov diffusion
processes. For Markov diffusion models the value process can be
represented as a space-transformation  of an asset price process
by the value function. The problem is to establish the
differentiability properties of the value function from the fact
that the value process satisfies the corresponding BSDE. The role
of the bridge between these equations is played by the statements
describing all invariant space-transformations of diffusion
processes, studied in Chitashvili and Mania (1996) and formulated
here in the appendix, in a suitable case adapted to financial
market models.  This approach enables us to prove that there
exists a solution (in a certain sense) of the Bellman equation and
that this solution is differentiable (in a generalized sense)
under mild assumptions on the model coefficients. Although, in our
case, the generalized derivative at $t$  and the second order
generalized derivatives at $x$ do not exist separately in general
(we prove an existence of a generalized $L$-operator), these
derivatives do not enter in the construction of  optimal
strategies which are explicitly given in terms of the first order
derivatives of the value function. It should be noted that the
theory of viscosity solutions is usually applied to such problems
(see, e.g., El Karoui et al (1997)), but differentiability of the
value function is in general beyond the reach of this method.

We assume that the dynamics of the asset price process is
determined by the following system of stochastic differential
equations
\begin{align}\label{um1}
dS_t=&\diag(S_t)(\mu(t,S_t,R_t)dt+\sigma^l(t,S_t,R_t)dW^l_t) \\
\label{um2} dR_t=&b(t,S_t,R_t)dt+\delta(t,S_t,R_t)dW^l_t+
\sigma^\bot(t,S_t,R_t)dW^\bot_t
\end{align}
Here $W=(W^1,..., W^n)$ be an $n$-dimensional standard Brownian
motion defined on a complete probability space $(\Omega,F, P)$
equipped with the $P$-augmentated filtration generated by $W$,
$F=(F_t,t\in[0,T])$. By $W^l=(W^1,...,W^d)$ and
$W^\bot=(W^{d+1},...,W^n)$ are denoted the $d$ and $n-d$
dimensional Brownian motions respectively.

 Assume that

{\bf S1)} the coefficients $\mu, b, \delta, \sigma^l, \sigma^\bot$
are measurable and bounded;

{\bf S2)} $n\times n-$ matrix function $\sigma\sigma'$ is
uniformly elliptic, i.e., there is a constant $c>0$ such that $$
(\sigma(t,s,r)\lambda,\sigma(t,s,r)\lambda)\ge c|\lambda|^2 $$ for
all $t\in[0,T],s\in R_+^d,r\in R^{n-d}$ and $\lambda\in R^n$,
where $\sigma$ is defined by
$$ \sigma(t,s,r)=
\left(\begin{array}{ll} \sigma^l(t,s,r) & 0 \\ \delta(t,s,r) &
\sigma^\bot(t,s,r),
\end{array}
\right).
$$
In addition we assume that

{\bf S3)} the system (\ref{um1}), (\ref{um2}) admits a unique
strong solution.

Straightforward calculations yield that in this case
$$\lambda=\diag(S)^{-1}(\sigma^l\sigma^{l'})^{-1}\mu,$$
where $\sigma^{l'}$ denotes the transpose of $\sigma^l$,
$$
\frac{d\la M\ra_t}{dt}=
diag(S_t)(\sigma^l\sigma^{l'})(t,S_t,R_t)diag(S_t)
$$
is the $\nu_t$ process, $\theta=(\sigma^l)^{-1}\mu$ is the market
price of risk and
$$
\la\lambda\cdot M\ra_t=\int_0^t||\theta_s||^2ds
$$
is the mean variance tradeoff.

By results of Krylov (1980) for sufficiently smooth coefficients
$\mu,\sigma,b,\delta$ the value function $V(t,x)$ can be
represented as $v(t,x,S_t,R_t)$  with sufficiently smooth function
$v(t,x,s,r),\; t\in [0,T],\; x\in R_+,\; s\in R_+^d,\;r\in
R^{n-d}$. Hence by the equation (\ref{ut20}) and the It\^o formula
we obtain that $v(t,x,s,r)$ satisfies the PDE

\begin{equation}\lbl{um1.1} {\mathcal
L}v(t,x,s,r)+v_s(t,x,s,r)'\diag(s)\mu(t,s,r)
+v_r(t,x,s,r)'b(t,s,r) \end{equation}

$$
=\frac{1}{2} \frac{|v_{sx}(t,x,s,r)
+diag(s)^{-1}\sigma^{l'}(t,s,r)^{-1}
\delta'(t,s,r)v_{rx}(t,x,s,r)+\lambda'(t,s,r)v_x(t,x,s,r)|_{\nu_t}^2}
{v_{xx}(t,x,s,r)},
$$
\begin{equation}\lbl{um1.2}
v(T,x,s,r)=U(x),
\end{equation}
which coincides with the Bellman equation of optimization problem
(\ref{ut3}), (\ref{um1}),(\ref{um2}) for controlled Markov
process. Moreover the optimal strategy is
$$
\pi^*(t,x,s,r)=
$$
$$
\frac{v_{sx}(t,x,s,r) +\diag(s)^{-1}\sigma^{l'}(t,s,r)^{-1}
\delta'(t,s,r)v_{rx}(t,x,s,r)+\lambda'(t,s,r)v_x(t,x,s,r)}
{v_{xx}(t,x,s,r)},
$$
In this section we study the solvability of
(\ref{um1.1}), (\ref{um1.2}) in the particular cases of utility
functions but with weaker conditions on coefficients.

First consider the case of power utility.

\be{th}\lbl{utt6.1} Let condition S1), S2) and S3) be satisfied.
Then the value function $v(t,s,r)$ admits all first order
generalized derivatives $v_s$ and $v_r$, a generalized L-operator
\begin{align*}
{\mathcal
L}v=&v_t+\frac{1}{2}tr(diag(s)\sigma^l\sigma^{l'}(t,s,r)\diag(s)v_{ss}+
tr(\delta\sigma^{l'}(t,r,s)diag(s)v_{sr}) \\
+&\frac{1}{2}tr((\delta\delta'(t,s,r)+
\sigma^\bot\sigma^{\bot'}(t,s,r))v_{rr})
\end{align*}
(in the sense of Definition D of the Appendix) and is the unique
bounded solution of equation
\begin{align}\label{bel}
\notag {\mathcal L}v&(t,s,r)+v_s(t,s,r)'diag(s)\mu(t,s,r)+
v_r(t,s,r)'b(t,s,r)\\
\notag =\frac{q}{2}&\frac{|v_{s}(t,s,r)
+diag(s)^{-1}\sigma^{l'}(t,s,r)^{-1} \delta'(t,s,r)v_{r}(t,s,r)
+\lambda(t,s,r)v(t,s,r)|_{\nu_t}^2}
{v(t,s,r)}\\
&dtdsdr-a.e.
\end{align}
with the boundary condition
\begin{equation}\label{bo}
v(T,s,r)=1.
\end{equation}
Moreover, the optimal strategy is defined as
$$
\pi^*(t,x,s,r)=(1-q)(\lambda(t,s,r)+\frac{\varphi(t,s,r)}{v(t,s,r)}\big)x
$$
and the optimal wealth process is of the form
$$
X_t^{*}=x{\mathcal E}_t((1-q)(\frac{\varphi}{v}+\lambda)\cdot S),
$$
 where
$\varphi(t,s,r)=v_{s}(t,s,r) +\diag(s)^{-1}\sigma^{l'}(t,s,r)^{-1}
\delta'(t,s,r)v_{r}(t,s,r)$. \ee{th}

{\it Proof.} {\it Existence.} Since $(S,R)$ is a Markov process,
the feedback controls are sufficient and the value process is
expressed by
\begin{equation}\label{um4}
V_t=v(t,S_t,R_t) \:\:\:  a.s.
\end{equation}where
$$
v(t,s,r)=\sup_{\pi\in\Pi_1}
E\big(\big(1+\int_t^T\pi_uds_u\big)^p|S_t=s,R_t=r\big).
$$
(one can show this fact, e.g., similarly to \cite{Ch5}).

\noindent Since the value process satisfies equation (\ref{power}),
it is an It\^o process. From the equality \\
${\mathcal E}(-\lambda\cdot M)={\mathcal E}(-\int_0^\cdot\theta_udw^l_u)$ and
boundedness of $\theta$ follows that ${\mathcal E}(-\lambda\cdot
M)$ satisfies The Muckenhoupt inequality. Thus
$V_t=\underset{\pi\in\Pi_1}{\text{\rm esssup}}\text{ }
E((1+\int_t^T\pi_udS_u )^p/{\mathcal F}_t) $ is bounded (see Remark 4.1)and
the
martingale part of $V$ is in BMO by Proposition 7 from \cite{M-sch}. Hence the
finite variation part of $V_t$ is of integrable variation and from
(\ref{um4}) we have that $v(t,S_t,R_t)$ is an It\^o process of the
form (\ref{D1}) (Appendix). Therefore, Proposition B of the
Appendix implies that the function $v(t,s,r)$ admits a generalized
L-operator, all first order generalized derivatives and can be
represented as
$$
v(t,S_t,R_t)=v_0+\int_0^t(v_s(u,S_u,R_u)'\diag(S_u)\sigma^l(u,S_u,R_u)
$$
$$
+v_r(u,S_u,R_u)'\delta(u,S_u,R_u))dW^l_s
+\int_0^tv_r(u,S_u,R_u)'\sigma^\bot(u,S_u,R_u)dW^\bot_s
$$
$$
+\int_0^t{\mathcal
L}v(u,S_u,R_u)ds+\int_0^t\big(v_s(u,S_u,R_u)'\diag(X_s)\mu(u,S_u,R_u)
$$
\begin{equation}\label{um5}
+ v_{r}(u,S_u,R_u)b(u,S_u,R_u)\big)du,
\end{equation}
where ${\mathcal L}V$ is the generalized ${\mathcal L}$-operator.

On the other side the value process is a solution of (\ref{power})
and by the uniqueness of the canonical decomposition of
semimartingales, comparing the martingale parts of (\ref{um5}) and
(\ref{power}), we have that $dt\times dP$- a.e.
\begin {align}\label{um6}
\varphi_t=&v_s(t,S_t,R_t)
+\diag(S_t)^{-1}\sigma^{l'}(t,S_t,R_t)^{-1}\delta'(t,S_t,R_t)v_r(t,S_t,R_t),\\
\label{m7}
\varphi^\bot_t=&\sigma^{\bot'}(t,S_t,R_t)v_r(t,S_t,R_t).
\end{align}
Then, equating the processes of bounded variation of the same
equations, taking into account (\ref{um5}) and (\ref{um6}), we
derive
\begin{align}
\notag\int_0^t
\big({\mathcal L}v(u,S_u,R_u)+v_s(u,S_u,R_u)'\diag(S_u)&\mu(u,S_u,R_u)\\
\notag
+v_{r}(u,S_u,R_u)b(u,S_u,R_u)\big)&du\\
=\frac{q}{2}\int_0^t\frac{|\varphi_u
+\lambda(u,S_u,R_u)v(u,S_u,R_u)|_{\nu_u}^2}{v(u,S_u,R_u)}
 &du
\end{align}
 which
gives that $v(t,s,r)$ solves the Bellman equation (\ref{bel}).

{\it Unicity.} Let $\tilde v(t,s,r)$ be a bounded positive
solution of (\ref{bel}), (\ref{bo}), from the class $V^L$. Then
using the generalized It\^o formula (Proposition B of Appendix)
and equation (\ref{bel}) we obtain that $\tilde v(t,S_t,R_t)$ is a
solution of (\ref{power}), hence $\tilde v(t,S_t,R_t)$ coincides
with the value process $v$ by Theorem 4.1. Therefore
$\tilde v(t,S_t,R_t)=v(t,S_t,R_t)$ a.s. and $\tilde v=v$, $dtdxdy$
a.e. \qed

Now we consider extreme cases for the stochastic volatility
models. Let first assume that coefficients
$\mu,\sigma^l$ does not contain the variable $r$. Hence  equation
(\ref{um1}) takes the form
\begin{equation}\label{mm1}
dS_t=\diag(S_t)(\mu(t,S_t)dt+\sigma^l(t,S_t)dW^l_t).
\end{equation}
Let $S(q)$ be the It\^o process governed by SDE
\begin{equation}\label{mm2}
dS_t(q)=\diag(S_t(q))\sigma^l(t,S_t(q))(dW^l_t+q\theta(t,S_t(q))dt),
\end{equation}
where $dW^l_t+q\theta(t,S_t)dt$ is Brownian motion w.r.t. measure
$dQ(q)={\mathcal E}_T(-q\int_0^\cdot\theta_u dw_u^l)dP.$ Thus  by
 \ref{4.6} the value process is represented as
$$
V_t=v(t,S_t(q))=(\tilde v(t,S_t(q))^{\frac{1}{1-q}},
$$
where $\tilde
v(t,s)=E^{Q(q)}\big(e^{\frac{q(q-1)}{2}\int_t^T|\theta_u|^2du}|S_t(q)=s\big).
$ Therefore we have

\be{cor}\lbl{utc6.1} Let conditions S1), S2) and S3) be satisfied
for the coefficients of system (\ref{mm2}). Then the value process
can be represented as $(\tilde v(t,S_t(q))^{\frac{1}{1-q}},$ where
$\tilde v(t,s)$ is the classical solution of the linear PDE
\begin{align}
\notag \tilde
v_t(t,s)+\frac{1}{2}tr(\diag(s)\sigma^l\sigma^{l'}(t,s)\diag(s)\tilde
v_{ss}(t&,s))
\\
+\frac{q(q-1)}{2}|\theta(t,s)|^2\tilde v(t,s)&=0,\\
\tilde v(T,s)&=1.
\end{align}
\ee{cor}

The second extreme case corresponds to the stochastic volatility
model of the form
\begin{align}
\notag
dS_t=&\diag(S_t)(\mu(t,S_t,R_t)dt+\sigma^l(t,S_t,R_t)dW^l_t) \\
\label{mex} dR_t=&b(t,R_t)dt+ \sigma^\bot(t,R_t)dW^\bot_t.
\end{align}

\be{cor}\lbl{utc6.2}  Let conditions S1), S2) and S3) be satisfied
for the coefficients of the system (\ref{mex}) and $\theta$ does not
depend on the variable $s$. Then the value process of the optimization
problem (\ref{pw}) is of the form $V_t=v(t,R_t)$, where
$v(t,r)=E(e^{-\frac{q}{2}\int_t^T|\theta(u,R_u)|^2du}|R_t=r)$
satisfies the linear PDE
\begin{align}
\notag
v_t(t,r)+\frac{1}{2}tr(\sigma^\bot\sigma^{\bot'}(t,r)v_{rr}(t,r))
\\
+v_r(t,r)'b(t,r)-\frac{q}{2}|\theta(t,r)|^2v(t,r)&=0,\\
v(T,r)&=1.
\end{align}
\ee{cor}

Similar results can be obtained for exponential utility
function.

\be{prop}\lbl{utp6.3} Let conditions S1), S2) and S3) be satisfied
and $H=g(S_T,R_T)$ for a continuous bounded function $g(s,r)$. Then
the value function $v(t,s,r)$ for the problem (\ref{exp}) admits
all first order generalized derivatives $v_s$ and $v_r$, a
generalized L-operator and is the unique bounded solution of
equation
\begin{align}\label{bell}
\notag {\mathcal L}v&(t,s,r)+v_s(t,s,r)'\diag(s)\mu(t,s,r)+
v_r(t,s,r)'b(t,s,r)\\
\notag =\frac{1}{2}&\frac{|v_{s}(t,s,r)
+diag(s)^{-1}\sigma^{l'}(t,s,r)^{-1} \delta'(t,s,r)v_{r}(t,s,r)
+\lambda(t,s,r)v(t,s,r)|_{\nu_t}^2}
{v(t,s,r)}\\
&dtdsdr-a.e.
\end{align}
with the boundary condition
\begin{equation}\label{bol}
v(T,s,r)=e^{-\gamma g(s,r)}.
\end{equation}
Moreover, the optimal strategy is defined as
$$
\pi^*(t,x,s,r)=\frac{1}{\gamma}(\lambda(t,s,r)+\frac{\varphi(t,s,r)}{v(t,s,r)}\big)x,
$$
 where
$\varphi(t,s,r)=v_{s}(t,s,r) +diag(s)^{-1}\sigma^{l'}(t,s,r)^{-1}
\delta'(t,s,r)v_{r}(t,s,r)$ and the optimal wealth process is
defined by (\ref{frexp}). \ee{prop}

The following assertion, for the case of  logarithmic utility, follows
immediately from
Theorem \ref{utt4.2} and the Feynmann-Kac formula:

\be{prop}\lbl{utp6.4}  Let condition S1), S2) and S3) be satisfied
and $U(x)=\log x$. Then the value function can be represented as
$v(t,S_t,R_t),$ where $v(t,s,r)$ is unique solution of linear PDE
\begin{align}
\notag {\mathcal L}v(t,s,r)+
v_s(t,s,r)'\diag(s)\mu(t,&s,r)\\
+v_r(t,s,r)'b(t,s,r)
+|\theta(t,s,r)|^2v(t,s,r)&=0,\\
v(T,s,r)&=1
\end{align}
and the optimal strategy is $\pi^*(t,x,s,r)=\lambda(t,s,r)x$.
\ee{prop}

\

{
\appendix

\section{}

\

Let us show that the family
\begin{equation}\label{uta1}
\Lambda_t^\pi=E(U(x+\int_0^T\pi_udS_u)|{{\mathcal F}}_t),
\;\;\pi\in\Pi_x(\tilde\pi,t,T)
\end{equation}
satisfies the $\vps$-lattice property (with $\vps=0$) for any
$t\in[0,T]$ and $\tilde\pi$. $\Pi(\tilde\pi,t,T)$ is a set of
predictable $S$-integrable processes $\pi$ from $\Pi_x$ such that
$$
\pi_s=\td\pi_sI_{(0\le s<t)}.
$$
We shall write $\Pi(t,T)$ instead of $\Pi(0,t,T)$ for the class of
strategies corresponding to $\tilde\pi=0$ up to time $t$.

We must show that for any $\pi^1,\;\pi^2\in\Pi(\td\pi,t,T)$ there
exists a strategy $\pi\in\Pi(\td\pi,t,T)$ such that
\begin{equation}\label{uta2}
\Lambda_t^\pi=\max(\Lambda_t^{\pi^1},\Lambda_t^{\pi^2}).
\end{equation}
For any $\pi^1$ and $\pi^2$ let us define the set
$$
B=\{\omega:\Lambda_t^{\pi^1}\le\Lambda_t^{\pi^2}\}
$$
and let
$$
\pi_s=\tilde\pi_sI_{(0\le s<t)}+ \pi_s^1I_BI_{(s\ge
t)}+\pi_s^2I_{B^c}I_{(s\ge t)}
$$
It is evident that
\begin{equation}\label{pipi}
\text{if}\;\;\;\tilde\pi, \pi^1,\pi^2\in\Pi_x,\;\;\text{then}\;\;\;
\pi\in\Pi_x.
\end{equation}

Since $B$ is ${\mathcal F}_t-$measurable we have
$$
\Lambda_t^\pi=E(U(x+\int_0^T\pi_udS_u)|{{\mathcal F}}_t)
$$
$$
=E(U(x+\int_0^t\td\pi_udS_u+I_B\int_t^T{\pi_u^1}dS_u+
I_{B^c}\int_t^T{\pi_u^2}dS_u)|{\mathcal F}_t)
$$
$$
=I_BE(U(x+\int_0^t\td\pi_udS_u+\int_t^T{\pi_u^1}dS_u)|{\mathcal
F}_t)+ I_{B^c}E(U(x+\int_\tau^t\td\pi_udS_u+
\int_t^T{\pi_u^2}dS_u)|{\mathcal F}_t)
$$
$$
=I_BE(U(x+\int_0^T{\pi_u^1}dS_u)|{\mathcal F}_t)+ I_{B^c}E(U(x+
\int_0^T{\pi_u^2}dS_u)|{\mathcal F}_t)
$$
$$
=E(U(x+\int_0^T{\pi_u^1}dS_u)|{\mathcal F}_t)\vee
E(U(x+\int_0^T{\pi_u^2}dS_u)|{\mathcal F}_t),
$$
hence (A.2) is satisfied.

{\bf Proposition A1)} (Optimality principle). Let condition B1) be
satisfied.

a) For all $x\in R$, $\pi\in\Pi$ and $s\in[0,T]$ the process
$(V(t,x+\int_s^t\pi_udS_u), t\ge s)$ is a supermartingale,
admitting an RCLL modification.

b) $\pi^*(s,x)$ is optimal iff $(V(t,x+\int_s^t{\pi^*_u}dS_u),
t\ge s)$ is a martingale.

c) for all $s<t$
\begin{equation}\label{uta3}
V(s,x)=\underset{\pi\in\Pi(s,T)}{\esssup}
E(V(t,x+\int_s^t\pi_udS_u)|{\mathcal F}_s).
\end{equation}

{\it Proof}. a) For simplicity we shall take $s$ equal to zero.
Let us show that $Y_t=V(t,x+\int_0^t\tilde\pi_udS_u)$ is
supermartingale for all $x$ and $\td\pi$. Since
$$
Y_t= \underset{\pi\in\Pi(t,T)}{\esssup}
E(U(x+\int_0^t\td\pi_udS_u+\int_t^T\pi_udS_u)|{\mathcal F}_t)
$$
using the lattice property of the family (A.1) from Lemma 16.A.5
of \cite{E} we have
$$
E(Y_t|{\mathcal F}_s)=E(\underset{\pi\in\Pi(t,T)}{\esssup}
E(U(x+\int_0^t\td\pi_udS_u+\int_t^T\pi_udS_u)|{\mathcal
F}_t)|{\mathcal F}_s)
$$
$$
=E(\underset{\pi\in\Pi(\td\pi,t,T)}{\esssup}
E(U(x+\int_0^T\pi_udS_u)|{\mathcal F}_t)|{\mathcal F}_s)
$$
\begin{equation}\label{uta4}
=\underset{\pi\in\Pi(\td\pi,t,T)}{\esssup}
E(U(x+\int_0^T\pi_udS_u)|{\mathcal F}_s).
\end{equation}
It is evident that $\Pi(\td\pi,t,T)\subseteq \Pi(\td\pi,s,T)$ for
$s\le t$, which implies the inequality
$$
\underset{\pi\in\Pi(\td\pi,t,T)}{\esssup}
E(U(x+\int_0^T\pi_udS_u)|{\mathcal F}_s)
$$
$$
\le \underset{\pi\in\Pi(\td\pi,s,T)}{\esssup}
E(U(x+\int_0^T\pi_udS_u)|{\mathcal F}_t)
$$
\begin{equation}\label{uta5}
=V(s,x+\int_0^s\tilde\pi_udS_u).
\end{equation}

Thus (A.4) and (A.5) imply that $E(Y_t/F_s)\le Y_s$.

b) If $V(t,x+\int_0^t\pi^*_udS_u))$ is a martingale, then

$$
\inf_{\pi\in\Pi}
EU(x+\int_0^T\pi_udS_u)=V(0,x)=EV(0,x)
$$
$$
=EV(T,x+\int_0^T{\pi^*_u}dS_u)
=EU(x+\int_0^T{\pi^*_u}dS_u),
$$
hence, $\pi^*$ is optimal.

Conversely, if $\pi^*$ is optimal, then
$$
EV(0,x)=\sup_{\pi\in\Pi}
EU(x+\int_0^T\pi_udS_u)
$$
$$
=EU(x+\int_0^T{\pi^*_u}dS_u)=EV(T,x+\int_0^T{\pi^*_u}dS_u).
$$
Since $V(t,x+\int_0^t{\pi^*_u}dS_u)$ is a supermartingale, the
latter equality implies that this process is a martingale (it
follows from Lemma 6.6 of \cite{L-Sh}).

c) Since $Y_t=V(t,x+\int_s^t\td\pi_udS_u)$ is a supermartingale
for any $\tilde\pi\in\Pi(s,T)$, $x\in R$ and $t\ge s$ we have
$$
V(s,x)\ge E(V(t,x+\int_s^t\td\pi_udS_u)|{\mathcal F}_s),
$$
hence
\begin{equation}\label{uta5}
V(s,x)\le \underset{\td\pi\in\Pi(s,T)}{\esssup}
E(V(t,x+\int_s^t\td\pi_udS_u)|{\mathcal F}_s).
\end{equation}
On the other hand for any $\td\pi$
$$
E(V(t,x+\int_s^t\td\pi_udS_u)|{\mathcal F}_s)=
$$
$$
E(\underset{\pi\in\Pi(t,T)}{\esssup}\;
E(U(x+\int_s^t\td\pi_udS_u+\int_t^T\pi_udS_u)|{\mathcal
F}_t){\mathcal F}_s)\ge
$$
$$
E(E(U(x+\int_s^T\td\pi_udS_u)|{\mathcal F}_t){\mathcal F}_s)=
E(U(x+\int_s^T\td\pi_udS_u)|{\mathcal F}_s).
$$
Taking esssup of the both parts we obtain
$$
\underset{\td\pi\in\Pi(s,T)}{\esssup}
E(V(t,x+\int_s^t\td\pi_udS_u)|{\mathcal F}_s)\ge
$$
\begin{equation}\label{uta7}
\underset{\td\pi\in\Pi(s,T)}{\esssup}
E(U(x+\int_s^T\td\pi_udS_u)|{\mathcal F}_s)=V(s,x).
\end{equation}
Thus the equality (A.3) follows from (A.6) and (A.7).

Let us show now that the process $\tilde
V(t,x+\int_0^t\td\pi_udS_u)$ admits an RCLL modification for each
$x\in R$ and $\pi\in\tilde\Pi$. According to Theorem 3.1 of
\cite{L-Sh} it is sufficient to prove that the function $E\tilde
V(t,x+\int_0^t\td\pi_udS_u)), t\in[0,T])$ is right-continuous for
every $x\in R$.

Let $(t_n,n\ge1)$ be a sequence of positive numbers such that
$t_n\downarrow t$, as $n\to\infty$. Since $\tilde
V(t,x+\int_0^t\td\pi_udS_u)$ is a supermartingale, we have
\begin{equation}\label{uta8}
 E\tilde V(t,x+\int_0^t\td\pi_udS_u)
\ge\lim_{n\to\infty}E\tilde V(t_n,x+\int_0^{t_n}\td\pi_udS_u).
\end{equation}
Let us show the inverse inequality. For $s=0$ equality (A.4) takes
the form
\begin{equation}\label{uta9}
E\tilde V(t,x+\int_0^t\td\pi_udS_u)=
\max_{\pi\in\tilde\Pi(\td\pi,t,T)} E(U(x+\int_0^T\pi_udS_u).
\end{equation}
Therefore, for any $\varepsilon>0$ there exists a strategy
$\pi^\varepsilon$ such that
\begin{equation}\label{uta10}
E\tilde V(t,x+\int_0^t\td\pi_udS_u)\le
E(U(x+\int_0^t\td\pi_udS_u+\int_t^T\pi^\varepsilon_udS_u)+\varepsilon.
\end{equation}
 Let us define a sequence $(\pi^n, n\ge1)$ of strategies
$$
\pi^n_s=\td\pi_sI_{(s<t_n)}+\pi_s^\varepsilon I_{(s\ge t_n)}.
$$
Using inequality (A.11), the continuity of $U$ (it follows from
B1) and B2), the convergence of the stochastic integrals and
Fatou's lemma, we have
$$
E\tilde V(t,x+\int_0^t\td\pi_udS_u)\le
E(U(x+\int_0^t\td\pi_udS_u+\int_t^T\pi^\varepsilon_udS_u)+\varepsilon=
$$
$$
=E(\lim_nU(x+\int_0^{t_n}\td\pi_udS_u+\int_{t_n}^T\pi^\varepsilon_udS_u))
+\varepsilon\ge
$$
$$
\ge\underline{\lim}_n E( E(U(x+\int_0^{t_n}\td\pi_udS_u+
\int_{t_n}^T\pi^\varepsilon_udS_u)/ {\mathcal
F}_{t_n}))+\varepsilon\ge
$$
$$
\ge\underline{\lim}_{n}
E(\underset{\pi\in\td\Pi(\td\pi,t_n,T)}{\esssup}
E(U(x+\int_0^{t_n}\td\pi_udS_u+\int_{t_n}^T\pi_udS_u)/{\mathcal
F}_{t_n}))+ \varepsilon=
$$
\begin{equation}\label{uta11}
=\underline{\lim}_{n\to\infty}E(\tilde V(t_n,
x+\int_0^{t_n}\td\pi_udS_u)+\varepsilon.
\end{equation}

Since $\varepsilon$ is an arbitrary positive number,  from (A.12)
we obtain that
\begin{equation}\label{uta12}
E\tilde
V(t,x+\int_0^t\td\pi_udS_u))\le\underline{\lim}_{n\to\infty}
E\tilde V(t_n,x+\int_0^{t_n}\td\pi_udS_u)),
\end{equation}
which together with (A.9) implies that the function $(E\tilde
V(t,x+\int_0^t\td\pi_udS_u)), t\in[0,T])$ is right-continuous.

\be{lem}\lbl{lemc}  Let $b_t$ be a predictable process and $S$ a
continuous semimartingale. Suppose that $\cal K$ is an adapted continuous increasing process and
$\mu^{\cal K}$  the corresponding  Dolean's measure. Denote by $\Pi_x$ the space of all
predictable $S$-integrable processes $\pi$ such that for all
$t\in[0,T]$
$$
x+\int_0^t\pi_udS_u\ge0.
$$
Then $\mu^{\cal K}\;\;a.e.$
$$
\underset{\pi\in\Pi_x}{\essinf} \big|\pi_t - b_t\big|=
0
$$
\ee{lem}
{\it Proof.}
Taking a bounded continuous approximation $b_t^{n,m}$ of
$b^n_t=b_tI_{(|b_t|\le n)}$ in the sense of $\mu^{\cal K}$-a.e. convergence we have that
$$
\underset{\pi\in\Pi_x}{\essinf} \big|\pi_t - b_t\big|\le
 \underset{\pi\in\Pi_x}{\essinf}
\big|\pi_t - b^{n,m}_t\big| + |b^{n,m}_t-b_t|.
$$
Therefore, without loss of generality we may assume that $b$ is
continuous and $S$-integrable.
Let us denote by
$$
\pi^{r,n}_t=b_t\psi_t^{r,n}{\mathcal E}_t(\frac{1}{x}(b\psi^{r,n})\cdot S),
$$
for $r\in[0,T],n\in N$, where
$\psi_t^{r,n}=I_{(r-\frac{1}{n},r+\frac{1}{n})}(t)(1-n|t-r|)$.
It is evident that $\pi^n$ belongs to $\Pi_x$ for all $r,n\ge1$.
Indeed,
$$
x+\int_0^t\pi_u^ndS_u=x + x\int_0^t{\mathcal E}_u (\frac{1}{x}b
\psi^{r,n}\cdot S)\frac{1}{x}b_u \psi^{r,n}(u)dS_u
$$
$$
= x{\mathcal E}_t(\frac{1}{x}b \psi^{r,n}\cdot S)\ge0.
$$
Denote by $\gamma_t$ the expression
$\underset{\pi\in\Pi_x}{\essinf} \big|\pi_t - b_t\big|$.
By definition  $\gamma_t\le |\pi_t-b_t|$, $\mu^{\cal K}$-a.e. for all $\pi\in\Pi_x$.
Therefore $\gamma_t\le |\pi^{r,n}_t-b_t|$ on the set $B$ with $\mu^{\cal K}(B^c)=0$ for all rational $r\in[0,T]$
and integer $n$.
Let
$$
\wt\gamma_t=\gamma_t,\;\;\text{if}\;\;(\omega,t)\in B\;\;\text{and}\;\;\wt\gamma_t=0,
\;\;\text{if}\;\;(\omega,t)\in B^c.
$$
Then $\wt\gamma_t\le |\pi^{r,n}_t-b_t|$ for all $t,\omega,r,n$. It
is easy to see that $\pi_t^{r,n}$ is continuous function of
variables $s,t$ for each $n$, since
$\psi_t^{s,n}\equiv\int_0^t(1_{(s-\frac{1}{n},s)}(u)-1_{(s,s+\frac{1}{n})}(u))du$
and
$$
\int_0^t\psi_u^{s,n}b_udS_u\equiv\psi_t^{s,n}\int_0^tb_udS_u-\int_0^t\left(\int_0^v
b_udS_u\right)\left(1_{(s-\frac{1}{n},s)}(v)-1_{(s,s+\frac{1}{n})}(v)\right)dv
$$
$$
\equiv\psi_t^{s,n}\int_0^tb_udS_u-\int_{(s-\frac{1}{n})\vee 0
}^{s\wedge t} b_udS_u+\int_{s\vee 0 }^{(s+\frac{1}{n})\wedge t}
b_udS_u
$$
are continuous.
Hence $\pi_t^{r,n}\to\pi_t^{s,n}$, as $r\to s$ uniformly in
$t\in[0,T]$. Passing to the limit as $r\to s$ we have
$$\wt\gamma_t\le |\pi^{s,n}_t-b_t|,\;\;\text{for all}\;\; t,s\;\;P-a.s.$$

Since $P-a.s.$
$$\pi_s^{s,n}=b_s{\cal E}_s(\frac{1}{x}(b\psi^{s,n})\cdot S)\rightarrow b_s,\;\;\;\;
\text{as}\;\;{n\to \infty},$$
we can conclude that $\wt\gamma_s=0$ for all $s$ $P-$a.s.. This implies that $\gamma_s=0$,
$\mu^{\cal K}-$a.e..

\

\section{}

Now we introduce some notions which enable us to present an
application of Theorem 1 to the Markov case.

Consider the system of stochastic differential equations
(\ref{um1}), (\ref{um2}) and assume that conditions S1) and S2) are
satisfied. Under these conditions there exists a unique weak
solution of (\ref{um1}), (\ref{um2}), which is a Markov process and
its transition probability function admits a density
$p(s,(x_0,y_0),t,(x,y))$ with respect to the Lebesgue measure. We
shall use the notation $p(t,x,y) = p(0,(x_0, y_0), t, (x,y))$ for
the fixed initial condition $S_0=x_0, R_0 =y_0$.

\noindent Introduce the measure $\mu$ on the space $([0,T]\times
R_+^d\times R^{n-d},{\mathcal B}([0,T]\times R_+^d\times
R^{n-d}))$ defined by
$$ \mu(dt,dx, dy)=p(t,x,y)dtdxdy. $$

\noindent Let $C^{1,2}$ be the class of functions $f$ continuously
differentiable at $t$  and twice differentiable at $x, y$ on
$[0,T]\times R_+^d\times R^{n-d}$. For functions $f\in C^{1,2}$
the $L$ operator is defined  as
\begin{align*}
Lf=&f_t+tr(\frac{1}{2}diag(x)\sigma^l\sigma^{l'}diag(x)f_{xx})+
tr(\delta\sigma^{l'}diag(x)f_{xy}) \\
&+tr(\frac{1}{2}(\delta\delta'+ \sigma^\bot\sigma^{\bot'})f_{yy})
\end{align*}
\noindent where  $f_t, f_{xx}, f_{xy}$ and $f_{yy}$ are partial
derivatives of the function $f$, for which we use the  matrix
notations.

{\bf Definition B.} We shall say that a function
$f=(f(t,x,y),t\ge0,x\in R_+^d, y\in R^{n-d})$ belongs to the class
$V^L_{\mu}$ if there exists a sequence of functions $(f^n,n\ge1)$
from $C^{1,2}$ and measurable $\mu$-integrable functions $f_{x_i}
(i\le d)$, $f_{y_j} (d<j\le n )$  and $(Lf)$ such that $$
E\sup_{s\le T}|f^n(s,S_s, R_s)-f(u,S_u,R_u)|\to0,\:\text{as}\:
n\to\infty, $$
$$ \iint_{[0,T]\times R^d_+\times R^{n-d}}
(f^n_{x_i}(s,x,y)-f_{x_i}(s,x,y))^2x_i^2\mu(ds,dx,dy)
\to0,\:\:i\le d, $$
$$ \iint_{[0,T]\times R^d_+\times
R^{n-d}} (f^n_{y_j}(s,x,y)-f_{y_j}(s,x,y))^2\mu(ds,dx,dy)
\to0,\:\: d<j\le n,
$$
$$ \iint_{[0,T]\times R^d_+\times
R^{n-d}}|Lf^n(s,x,y)- (Lf)(s,x,y)|\mu(ds,dx,dy)\to0,
$$
$$
 \text{as}\;\; n\to\infty.
$$

Now we formulate the statement proved in Chitashvili and Mania
(1996) in the case convenient for our purposes.

{\bf Proposition B.} Let conditions S1)-S2) be satisfied and let
$f(t, S_t, R_t)$ be a bounded process. Then the  process
$(f(t,S_t,R_t), t\in[0,T])$ is  an It\^o process of the form $$
f(t,S_t,R_t)=f(0,S_0,R_0)+\int_0^tg(s,\omega)dW_s+
\int_0^ta(s,\omega)ds,\:\:\: {\text a.s.} $$ with
\begin{equation}\label{D1}
E\int_0^tg^2(s,\omega)ds<\infty,\:\:
E\int_0^t|a(s,\omega)|ds<\infty
\end{equation}
if and only if $f$ belongs to $V^L_{\mu}$. Moreover the process
$f(t,S_t, R_t)$ admits the decomposition $$ f(t,S_t, R_t)=
f(0,S_0,R_0)+\sum_{i=1}^d\int_0^tf_{x_i}(s,S_s,R_s)dS^i_s + $$
\begin{equation}
\sum_{j=d+1}^n\int_0^tf_{y_j}(s,S_s,R_s)dR^j_s +
\int^t_0(Lf)(s,S_s,R_s)ds.
\end{equation}

{\bf Remark.} For continuous functions $f \in V^L_{\mu}$ the
condition
\begin{equation}
\sup_{(t,x,y)\in D}|f^n(t,x,y)-f(t,x,y)|\to0,\text{ as }
n\to\infty
\end{equation}
for every compact $D\in [0,T]\times R^d_+\times R^{n-d}$, can be
used instead of the first relation  of Definition B.
}

\

\end{document}